\documentclass[leqno]{siamltex}
\usepackage{amsmath,mathtools}
\usepackage{amssymb}
\usepackage{graphicx}
\usepackage{color}
\usepackage{epstopdf} 
\usepackage{multirow}
\usepackage{tabularx}
\usepackage[ruled,vlined]{algorithm2e}

\newtheorem{remark}{Remark}

\newcommand{\cO}{\mathcal{O}}

\renewcommand{\div}{\mbox{\rm div\,}}

\newcommand{\mE}{\mathbb{E}}

\newcommand{\mP}{\mathbb{P}}

\newcommand{\A}{\mathbf{A}}
\newcommand{\B}{\mathbf{B}}

\newcommand{\U}{\mathbf{U}}
\newcommand{\bb}{\mathbf{b}}

\newcommand{\bu}{\mathbf{u}}
\newcommand{\bbf}{\mathbf{f}}

\newcommand{\Ome}{\Omega}
\newcommand{\p}{\partial}
\newcommand{\nab}{\nabla}

\begin{document}

\title{An efficient iterative method for solving parameter-dependent and random   convection-diffusion problems}
\markboth{X. FENG, Y. LUO, L. VO, and Z. WANG}{NUMERICAL METHOD FOR PARAMETER-DEPENDENT PROBLEMS}

\author{
Xiaobing Feng\thanks{Department of Mathematics, The University of Tennessee,
Knoxville, TN 37996, U.S.A. ({\tt xfeng@math.utk.edu})
The work of this author was partially supported by
the NSF grants DMS-1620168 and DMS-2012414.}
\and
Yan Luo\thanks{School of Mathematical Sciences, University of Electronic Science and Technology of China, Chengdu, Sichuan 611731, China ({\tt luoyan\_16@126.com}) The work of this author was supported by the Young Scientists Fund of the National Natural Science Foundation of China grant 11901078 and by the Fundamental Research Funds for the Central Universities grant  ZGX2020J021.}
\and
Liet Vo\thanks{Department of Mathematics, The University of Tennessee, Knoxville, TN 37996, U.S.A. ({\tt lvo6@vols.utk.edu} The work of this author was partially supported by
the NSF grants DMS-1620168 and DMS-2012414.}
\and
Zhu Wang\thanks{Department of Mathematics, University of South Carolina, Columbia, SC 29208, U.S.A. ({\tt wangzhu@math.sc.edu}) 
The work of this author was partially supported by the NSF grant DMS-1913073.}	
}

\maketitle

\begin{abstract}
This paper develops and analyzes a general iterative framework for solving parameter-dependent and 
random   convection-diffusion  problems. It is inspired by the multi-modes method of \cite{Feng,Feng2} and the ensemble method of \cite{Yan} and extends those methods into a more general and unified framework. The main idea of the framework is to reformulate the underlying problem into another problem with parameter-independent   convection and diffusion coefficients and a parameter-dependent (and solution-dependent) right-hand side, a fixed-point iteration is then employed to compute the solution of the reformulated problem. 
The main benefit of the proposed approach is that an efficient direct solver and a block Krylov subspace iterative solver can be used at each iteration, allowing to reuse the $LU$ matrix factorization or to do an efficient matrix-matrix multiplication for all parameters, which in turn results in  significant computation saving. Convergence and rates of convergence are established for the iterative method both at the variational continuous level and at the finite element discrete level under some structure conditions. Several strategies for establishing  reformulations of parameter-dependent and random diffusion and   convection-diffusion  problems are proposed and their computational complexity is analyzed.  Several 1-D and 2-D numerical experiments are also provided to demonstrate the efficiency of the proposed iterative method and to validate the theoretical convergence results.

\end{abstract}



\section{Introduction} \label{sec-1}

Parameter-dependent problems arise from many engineering and scientific applications such as 
fluid mechanics, porous media flow, wave propagation and various biological models 
(cf. \cite{WCR} and the references therein). 
The parameter could be deterministic (such as viscosity, diffusivity, permeability, frequency, and birth/death rate) or random,  they appear as 
parameters in their respective mathematical (i.e., PDE) models. As the solution of such a model 
depends on the parameter, usually in a nonlinear manner even the PDE is linear, to compute 
approximate solutions,  one is forced 
to solve the same problem multiple times for a given set of parameters. Such a direct (or brute force) approach
may still be feasible if the parameter set is small even it is not efficient,  however, it becomes too expensive to use if the parameter set is large due to the sheer amount of computation required 
to solve a complicated PDE problem tens, hundreds even thousands times. 

To overcome the computational challenge, several approaches have been proposed in literature including variants of the model order reduction (MOR) method \cite{hesthaven2015certified,quarteroni2015reduced,antoulas2020interpolatory,brunton2019data}, the ensemble method \cite{jiang2014algorithm,Gunzburger2016efficient,gunzburger2017ensemble,gunzburger2017second,Yan,luo2018multilevel} and the multi-modes method \cite{Feng,Feng2,Feng3}. The MOR methods aim to lower the computational complexity of the underlying model/problem by approximating the solution manifold using a handful of degrees of freedom, a low-dimensional surrogate to the original model is then built which is usually cheaper to simulate and often referred to as a reduced order model. Therefore, the computational 
saving is due to the dimensionality reduction. On the other hand, noticing the similarity of solving the same type problems with different parameters, a natural idea is to solve those problems together as a group and to reuse the computation as much as possible, so the computational saving is due to reusing a major portion of the computation for a set of parameters and requires designing better algorithms. This is exactly the idea used by the ensemble and multi-modes methods, and is also the approach adopted in this paper. 

Inspired by the ensemble and multi-modes methods, the primary objective of this paper is 
to develop and analyze a more general and unified iterative framework for solving parameter-dependent (including the case of random parameters)  convection-diffusion  problems in an 
abstract variational setting. Our main idea is to reformulate a given parameter-dependent problem into another problem with parameter-independent   convection  and diffusion coefficients and 
a parameter-dependent (and solution-dependent) right-hand side, its solution is
then computed using a fixed-point iteration algorithm. Since the   convection and  diffusion coefficients 
are parameter-independent and iterate-independent, while the right-hand side term is parameter-dependent and iterate-dependent, this presents an ideal set-up for us to employ an efficient direct solver and a block Krylov subspace iterative solver at each iteration for solving a linear system with multiple right-hand vectors. It is the reuse of the $LU$ matrix factorization and efficiently to compute matrix-matrix multiplication for all parameters that results in significant computation saving, hence, leads to an efficient overall numerical method.

The remainder of this paper is organized as follows. In Section \ref{sec-2} we present 
the abstract variational setting for the parameter-dependent  convection-diffusion problems to be 
considered in this paper. Structure conditions are stated to ensure the well-posedness of 
the variational problem. Such an abstract setting allows the wider applicability of the method and results of this paper, it also permits a clean and concise presentation of the convergence analysis. 
In Section \ref{sec-3} we first present our main algorithm in the variational continuous  
level and then establish the rate of convergence of the algorithm under some structure 
conditions. Section \ref{sec-4} presents finite element approximations of both the
variational problem and the iterative algorithm. The highlight of this section is 
to derive the rate of convergence for the discretized algorithm. It worths noting that 
the finite element approximations are chosen for the preciseness and ease of 
presentation, the convergence results are also valid for other Galerkin-type 
approximations such as discontinuous Galerkin and spectral approximations. 
Section \ref{sec-5} is devoted to computational complexity analysis of 
the proposed algorithm and to discussing an implementation strategy 
of the algorithm based on a CVT (centroidal Voronoi tessellations) clustering method.  Finally, in Section \ref{sec-6} we present  three  sets of numerical tests, the first one solves a  
parameter-dependent diffusion problem and the second solves a random diffusion problem,
  and the third computes a random convection-dominated convection-diffusion problem which is a randomized double-glazing problem \cite{IFISS}. 
Both 1-D and 2-D numerical test results are presented to demonstrate the efficiency 
of the proposed iterative method and to validate the theoretical convergence results.
In particular,   the last two tests show that the proposed iterative method 
can successfully solve some strongly random diffusion and  convection-diffusion problems which the multi-modes method fails to do.

\section{Variational problem setting}\label{sec-2}
Let $H$ be a Hilbert space and $V \subset H$ be a Banach space. $V^*$ denotes
the dual space of $V$, the space of bounded linear functionals on $V$.  We consider the following variational problem of seeking $u \in V$ such that 
\begin{align}\label{eq1.1}
	A(\omega; u,v) &= F(\omega;v) \qquad\forall v \in V,
\end{align}
where $F(\omega;\cdot)\in V^*$ and $A(\omega;\cdot,\cdot): V\times V \rightarrow \mathbb{R}$ depends on a parameter $\omega$, which belongs to a (given) parameter space $\Ome$ and satisfies the following properties:
\begin{itemize}
	\item[(i)] \textit{Bi-linearity:} $A$ is bilinear mapping on $V\times V$.
	\item[(ii)] \textit{Boundedness:} There exists an $\omega$-independent constant $\Lambda > 0$ such that $$|A(\omega; v,w)| \leq \Lambda \|v\|_V\|w\|_V \qquad \forall v,w \in V. $$
	\item[(iii)] \textit{Coercivity:} There is an $\omega$-independent  constant $\lambda>0$ such that $$A(\omega; v,v) \geq \lambda \|v\|^2_V \qquad \forall v \in V.$$
\end{itemize}
  
\smallskip
\begin{remark}
	\begin{itemize}
		\item[(a)]  Since $A(\omega,\cdot,\cdot)$ may be non-symmetric, then 
			problem \eqref{eq1.1} contains both coercive diffusion and coercive convection-diffusion problems which
			could be convection-dominated, see Section \ref{sec-6.3}-\ref{sec-6.4} for more discussions. 
		\item[(b)] The parameter space $\Ome$ can be a deterministic or probability space, in the latter case, equation \eqref{eq1.1} becomes a random PDE. Moreover, $\Ome$ can be finite or infinite set.  
		\item[(c)] Applicable examples to be considered later are (i) parameter-dependent PDEs; (ii) random PDE.
	\end{itemize}
\end{remark}

It is well known \cite{BrennerScott} that under the assumptions (i)--(iii), there holds the following 
Lax-Milgram Theorem.

\begin{proposition}\label{lem1.1}
	Under the above assumptions, problem  \eqref{eq1.1} has a unique solution $u(\omega)\in V$ 
	for each $\omega\in \Omega$ and there holds
	\begin{align}\label{stability}
		\|u(\omega)\|_V \leq \frac{\|F(\omega)\|_{V^*}}{\lambda}   \qquad\forall \omega\in \Omega.
	\end{align}
\end{proposition}
 
\begin{remark}
	 We note that $\|F(\omega)\|_{V^*}<\infty$ for each $\omega\in \Omega$, however, the bound may 
	 depend on $\omega$. In the subsequent sections, unless stated otherwise, we shall assume that there exists an $\omega$-independent constant $C_F >0$ such that $\|F(\omega)\|_{V^*}\leq C_F$
	 for all $\omega\in \Ome$,  that is, the operator norm of $F(\omega)$ is uniformly bounded in $\Ome$.  This condition will be
	 relaxed to $E[\|F(\cdot)\|_{V^*}] <\infty$ for random diffusion  and  convection-diffusion applications in Section \ref{sec-6}, namely, the expected value of the operator norm of $F(\omega)$ is bounded.
\end{remark}	

As mentioned earlier, the aim of this paper is to develop efficient numerical methods/algorithms for computing 
the solution set $u(\omega)$ for all $\omega\in \Ome$. 

\section{The proposed iterative framework}\label{sec-3}
The goal of this section is to formulate our main iterative algorithm for solving problem 
\eqref{eq1.1} and to establish the rate of convergence for the algorithm. 

\subsection{Formulation of the main algorithm} Let $A_0(\cdot,\cdot): V\times V \rightarrow \mathbb{R}$ be an  $\omega$-independent bilinear form, and set $A_1 = A - A_0$. Trivially
\begin{align}
	A(\omega;\cdot,\cdot) = A_0(\cdot,\cdot) + A_1(\omega;\cdot,\cdot) 
\end{align}
is also a bilinear form on $V\times V$. The equation \eqref{eq1.1} can be rewritten as 
\begin{align}
\label{eq2.2} A_0(u,v) = F(\omega; v) - A_1(\omega;u,v).
\end{align}

We now are ready to state our main iterative algorithm for solving \eqref{eq1.1}.

\begin{algorithm}[H]
\small
    \caption{Main Algorithm in Variational Setting}
    \label{alg:algorithm1}
    \DontPrintSemicolon

	Step 1: Find $U_0\in V$ by solving the following problem:
	\begin{align}\label{eq2.6}	
	A_0\bigl(U_0,v\bigr) &= F(\omega;v)\qquad\forall v \in V.
	\end{align}

Step 2: For each $\omega\in \Omega$, determine $\{U_n=U_n(\omega)\}_{n\geq 1}\subset V$ recursively by solving 
\begin{align}\label{eq2.7} 
A_0\bigl(U_n,v\bigr) = F(\omega;v) - A_1\bigl(\omega; U_{n-1},  v\bigr) \qquad\forall v \in V.
\end{align}
	
\end{algorithm}

To complete the construction of Algorithm 1, we need to specify the bilinear form $A_0$. As expected, it  must be problem-dependent, hence, we will do so in Section \ref{sec-6} when we apply Algorithm 1 to specific application problems. We conclude this subsection with the following
remarks.

\begin{remark}
	The main advantage of Algorithm 1 is that the left hand-side bilinear forms in \eqref{eq2.6} and \eqref{eq2.7} are the same and independent of the parameter $\omega$ which 
	allows to design fast solvers for computing the solutions in both steps of
	Algorithm \ref{alg:algorithm1}. 
\end{remark}

\begin{remark}
(a) To recover the multi-modes method of \cite{Feng} for random diffusion problems 
with the diffusion coefficients of the form $a(\omega, x)=a_0(x)+\eta(\omega, x)$
(see Section \ref{sec-6.2} for the details), let 
$a_0(x)$ be the expected value of $a(\omega, x)$ and we further assume that $\eta=\varepsilon \xi$, then 
define the $n$th mode function by 
\begin{equation}\label{u_n_mode}
u_0=U_0, \qquad u_n:= \frac{1}{\varepsilon^n} \bigl( U_n- U_{n-1} \bigr) \quad\mbox{for } n\geq 1,
\end{equation}
which implies that 
\[
U_n= U_{n-1} + \varepsilon^n u_n = u_0+\varepsilon u_1 + \varepsilon^2 u_2 \cdots + \varepsilon^n u_n.
\]
Thus, the multi-modes method of \cite{Feng} is recovered. We note that the multi-modes method of \cite{Feng} 
does compute the mode functions $\{u_j\}_{j=1}^n$ in order to form the final approximate $U_n$. 
On the other hand, Algorithm 1 does not compute the mode functions, instead, it computes 
the approximate $U_n$ directly.  We also note that the parameter space $\Omega$ does not 
need to be given a priori but can be sampled on fly in simulations. 

(b) To recover the ensemble method of \cite{Yan} for parameter-dependent problems with the
leading-term coefficients $a(\omega_j, x)$ and  
the parameter set $\Omega=\{\omega_j\}_{j=1}^J$ which are a priori given (see Section \ref{sec-6.1} for the details), we simply choose $a_0$ to be 
\[
a_0(x) = \frac{1}{J} \sum_{j=1}^J a(\omega_j, x) \quad\mbox{or }\quad
a_0(x) = \max_{1\leq j\leq J} a(\omega_j, x) 	
\qquad \forall x\in D,
\]
Algorithm \ref{alg:algorithm1} then leads to the ensemble method. 
\end{remark}

\subsection{Convergence analysis} \label{sec-3.2}
To ensure the convergence of Algorithm \ref{alg:algorithm1}, we impose the following criterion (and practical 
guideline) for choosing $A_0$ which would ensure the convergence of the algorithm.

\medskip
{\bf Criterion for selecting $A_0$}
\smallskip

\begin{itemize}
	\item[(a)] $A_0$ must satisfy (i), (ii), (iii) from definition of $A$. That is,
	\begin{itemize}
		\item[(i)] $A_0$ is bilinear on $V\times V$.
		\item[(ii)] There exists a constant $\Lambda_0>0$ such that $$|A_0(v,w)| \leq \Lambda_0\|v\|_V\|w\|_V
		\qquad \forall v,w \in V.$$
		\item[(iii)] There is a constant $\lambda_0>0$ such that $$ A_0(v,v) \geq \lambda_0 \|v\|_V^2
		\qquad \forall v \in V. $$
	\end{itemize}

	\item[(b)] \textit{Relative dominance}: There exists a $\rho \in (0,1)$ such that 
	\begin{align}
	\frac{\|A_1(\omega)\|_{\mathcal{L}(V\times V, \mathbb{R})}}{\lambda_0} < \rho\qquad \forall 
	\omega\in \Omega.
	\end{align}
\end{itemize}

Under the above assumptions, we have the following convergence theorem. 

\begin{theorem}\label{thm3.1}
	Let $u(\omega)$ be the solution of \eqref{eq1.1} and $\{U_n(\omega)\}_{n\geq 0}$ be generated by Algorithm \ref{alg:algorithm1}. Suppose the conditions of the above criterion hold. Then there exists 
	a constant $C>0$   independent of $\Ome$ such that 
\begin{equation} \label{convergence_rate}
\|u(\omega)-U_n(\omega)\|_V\leq C \rho^{n+1} \qquad \forall \omega\in \Ome.
\end{equation}
Hence, $U_n(\omega)$ converges to  $u(\omega)$ strongly in $V$ as  $n\rightarrow \infty$ for every $\omega\in \Omega$.
\end{theorem}

\begin{proof} 
	For each fixed $\omega\in \Omega$, let $r_n=r_n(\omega) = u(\omega) - U_n(\omega) \in V$ for $n \geq 0$. Then, subtracting equation \eqref{eq2.6}  from \eqref{eq2.2} yields the following error equation for $n=0$:
	\begin{align}
\label{eq2.14}		A_0\bigl(r_0, v\bigr) = -A_1(\omega; u,v) \qquad\forall v \in V.
	\end{align}
	
Similarly,  for any $n \geq 1$ subtracting \eqref{eq2.7} from \eqref{eq2.2} gives the error equation 
\begin{align}
\label{eq2.15}	A_0\bigl(r_n, v\bigr) = -A_1\bigl(\omega; r_{n-1}, v\bigr)\qquad\forall v \in V.
\end{align}
Choose $v = r_0$ in \eqref{eq2.14} we obtain:
\begin{align}
	\lambda_0\|r_0\|^2_V \leq A_0\bigl(r_0, r_0\bigr) &= -A_1\bigl(\omega; u, r_0\bigr) 
	\leq \|A_1(\omega)\|_{\mathcal{L}(V\times V,\mathbb{R})}\, \|u\|_V\|r_0\|_V.
\end{align}
This implies that
\begin{align}
	\|r_0\|_V \leq \frac{\|A_1(\omega)\|_{\mathcal{L}(V\times V, \mathbb{R})}}{\lambda_0}\, \|u\|_V \leq \rho \,\|u\|_V.
\end{align}

We also choose $v = r_n$ in \eqref{eq2.15} to obtain:
\begin{align}
	\lambda_0\|r_n\|_V^2 \leq A_0\bigl(r_n,r_n\bigr) &= - A_1\bigl(\omega;r_{n-1}, r_n\bigr) \\\nonumber
	&\leq \|A_1(\omega)\|_{\mathcal{L}(V\times V,\mathbb{R})}\,\|r_{n-1}\|_V\|r_n\|_V.
\end{align}
Therefore,
\begin{align}
	\|r_n\|_V \leq \frac{\|A_1(\omega)\|_{\mathcal{L}(V\times V, \mathbb{R})}}{\lambda_0}\,\|r_{n-1}\|_V \leq \rho \|r_{n-1}\|_{V}.
\end{align}
By induction we have
\begin{align}
	\|r_n\|_{V} \leq \rho^{n+1}\,\|u\|_{V} \leq \rho^{n+1} \,\frac{1}{\lambda}\,\|F\|_{V^*}. 
\end{align}
The proof is complete.
\end{proof}

\section{Finite element approximation} \label{sec-4}
In this section, we first formulate the finite element Galerkin approximation of  \eqref{eq1.1} 
and then present a discrete analogue of Algorithm~\ref{alg:algorithm1} as a fast solver for computing the finite element solutions for all $\omega\in \Omega$. To the end,  let $\mathcal{T}_h$ be a quasi-uniform triangulation of a bounded domain $D\subset \mathbb{R}^d\, (d=1,2,3)$ with mesh size $h \in (0,1)$. Let $V^h_r \, (r\geq 1)$ denote the finite element space consisting of continuous piecewise $r$th order polynomials associated with $\mathcal{T}_h$. We assume that $V^h_r$ is a subspace of $V$,
which implicitly assumes that $V\subset W^{1,1}(D)$ but $V\not\subset W^{2,1}(D)$. 

\smallskip

For each $\omega\in \Omega$, the standard finite element Galerkin approximation of \eqref{eq1.1} is defined as seeking $u_h=u_h(\omega) \in V^h_r$ such that 
\begin{align}\label{eq4.1}	
A\bigl(\omega; u_h,v_h\bigr) = F(\omega;v_h)\qquad\forall v_h\in V^h_r.
\end{align}
It is easy to show that \eqref{eq4.1} has a unique solution $u_h$ which 
also satisfies the stability estimate \eqref{stability}. Moreover, there holds
the following error estimate.
 
\begin{proposition}\label{lem4.1}
	Let $u(\omega)$ and $u_h(\omega)$ denote the solutions of \eqref{eq1.1} and \eqref{eq4.1} respectively. Then, there exists a constant $C=C(\omega)>0$ independent of $h$ and a positive 
	integer $\ell (\leq r)$ such that 
	\begin{align}
		\|u(\omega) - u_h(\omega)\|_{V} \leq C \,h^\ell \qquad \forall \omega\in \Ome.
	\end{align}
\end{proposition}

\begin{remark}
	The dependence of $C$ on $\omega$ is through the norm $\|u(\omega)\|_{H^{r+1}}$ and 
	$\ell$ depends on	the space $V$. For example, if $V=H^1_0(D)$, then $\ell=r$. 
\end{remark}

Using the definition of $A_0$ and $A_1$, we can rewrite \eqref{eq4.1} as follows:
\begin{align}
\label{eq4.2}	A_0\bigl(u_h,v_h\bigr) = F(\omega;v_h) - A_1\bigl(\omega; u_h,v_h\bigr) \qquad\forall v_h\in V_r^h.
\end{align}
We then can easily formulate the discrete counterpart of Algorithm \ref{alg:algorithm1} for 
computing the solution of \eqref{eq4.1} below.

\begin{algorithm}[H]
\small
    \caption{Main Algorithm in Discrete Setting}
    \label{alg:algorithm2}
    \DontPrintSemicolon

	Step 1: Find $U_0^h \in V_r^h$ by solving the following problem:
	\begin{align}\label{eq4.3}	
	A_0\bigl(U_0^h,v_h\bigr) &= F(\omega;v_h)\qquad\forall v_h \in V_r^h.	
	\end{align}

	Step 2: For each $\omega \in \Ome$, determine $\{U^h_n=U^h_n(\omega)\}_{n\geq 1}\subset V_r^h$ recursively by solving  
	\begin{align}\label{eq4.4}	
	A_0\bigl(U_n^h,v_h\bigr) = F(\omega;v_h) - A_1\bigl(\omega; U^h_{n-1},  v_h\bigr) \qquad\forall v_h \in V_r^h. 
	\end{align}
	
\end{algorithm}

\begin{remark}
	We note that the main advantage of Algorithm \ref{alg:algorithm2} is that the left hand-side bilinear forms (or stiffness matrices) 
	in \eqref{eq4.3} and \eqref{eq4.4} are the same and independent of the parameter $\omega$ which 
	allows to design fast solvers for computing the solutions in both steps of
	Algorithm \ref{alg:algorithm2}. 
\end{remark}

Similarly, there also holds the following convergence result for Algorithm \ref{alg:algorithm2}. 

\begin{theorem}\label{thm4.2}
	Let $u_h(\omega)$ be solution of \eqref{eq4.2} and $\{U_n^h(\omega)\}_{n\geq 0}$ be generated by Algorithm \ref{alg:algorithm2}. Then there exists a constant $C>0$ independent of $\omega$ such that
	\begin{align}
		\|u_h(\omega) - U_n^h(\omega)\|_V \leq C\rho^{n+1} \qquad\forall \omega\in \Ome.
	\end{align}
Hence, $U_n^h(\omega)$ converges to  $u^h(\omega)$ strongly in $V$ as  $n\rightarrow \infty$ for every $\omega\in \Omega$.
\end{theorem}

\begin{proof} 
For each fixed $\omega\in \Ome$, let $e_n^h = u_h(\omega) - U_n^h(\omega) \in V^h_r$ for $n \geq 0$. Then, subtracting equation \eqref{eq4.3} from \eqref{eq4.2} yields the following error equation for $n=0$:
\begin{align}
\label{eq4.7}	A_0\bigl(e_0^h,v_h\bigr) = -A_1\bigl(\omega;u_h,v_h\bigr)\qquad\forall v_h \in V_r^h.
\end{align}
Similarly, by subtracting \eqref{eq4.4} from \eqref{eq4.2}, we obtain the following error equation for $n\geq 1$:
\begin{align}\label{eq4.8}	
A_0\bigl(e_n^h, v_h\bigr) = -A_1\bigl(\omega; e_{n-1}^h, v_h\bigr) \qquad\forall v_h \in V_r^h.
\end{align}

Choosing  $v_h = e^h_0\in V_r^h$ in \eqref{eq4.7}, we have
\begin{align}
	\lambda_0\|e^h_0\|_V^2 \leq A_0\bigl(e_0^h,e_0^h\bigr) &= -A_1\bigl(\omega;u_h,e^h_0\bigr) 
	\leq \|A_1(\omega)\|_{\mathcal{L}(V\times V, \mathbb{R})} \|u_h\|_V\|e^h_0\|_V.
\end{align}
This implies that
\begin{align}
	\|e_0^h\|_V \leq \frac{\|A_1(\omega)\|_{\mathcal{L}(V\times V, \mathbb{R})}}{\lambda_0}\, \|u_h\|_V \leq \rho \|u_h\|_V.
\end{align}
In addition, setting $v_h = e_n^h \in V_r^h$ in \eqref{eq4.8} yields
\begin{align}
	\lambda_0\|e_n^h\|_V^2 \leq A_0\bigl(e_n^h,e_n^h\bigr) &= -A_1\bigl(\omega;e_{n-1}^h,e_n^h\bigr)\\\nonumber
	& \leq \|A_1(\omega)\|_{\mathcal{L}(V\times V,\mathbb{R})}\|e_{n-1}^h\|_V\|e_n^h\|_V.
\end{align}
Therefore,
\begin{align}
	\|e_n^h\|_V \leq \frac{\|A_1(\omega)\|_{\mathcal{L}(V\times V,\mathbb{R})}}{\lambda_0}\|e_{n-1}^h\|_V \leq \rho \|e_{n-1}^h\|_V.
\end{align}
By induction we have 
\begin{align}
	\|e_n^h\|_V \leq \rho^{n+1}\|u_h\|_V \leq \rho^{n+1} \frac{1}{\lambda}\|F\|_{V^*}.
\end{align}
The proof is complete.
\end{proof}

An immediate corollary is the following global error estimate.

\begin{theorem}\label{thm4.3}
	Let $u(\omega)$ be solution of \eqref{eq1.1} and $\{U_n^h(\omega)\}_{n\geq 0}$ be generated by Algorithm 2.  Then, there holds
	\begin{align}
		\|u (\omega)- U^h_n(\omega) \|_V \leq C\bigl(h^\ell + \rho^{n+1} \bigr) \qquad \forall \omega\in \Ome.
	\end{align}
\end{theorem}

\begin{proof}
	The proof follows from Proposition \ref{lem4.1}, Theorem \ref{thm4.2} and an application of  the triangular inequality to the decomposition $u - U_n^h = (u - u_h) + (u_h - U^h_n)$.
\end{proof}

\section{Computational complexity and implementation strategies}\label{sec-5} 

\subsection{Linear solvers}
Steps in Algorithm~\ref{alg:algorithm2} solve linear systems with a coefficient matrix independent of the parameter $\omega$, which is a highly appealing feature when a group of problems, associated to different $\omega$'s, needs to be solved. 
Since their discrete systems would share a common coefficient matrix, one can solve them simultaneously from a single system with multiple right-hand-side (RHS) vectors. 
Suppose $J$ problems are considered corresponding to $\{\omega_1,\ldots, \omega_J\}$. For $j$th problem, at the algebraic level, Algorithm~\ref{alg:algorithm2} finds $\bu_n^j$ such that: for $n\geq 0$, 
\begin{equation}
\A_0 \bu_n^j = \bb_n^j 
=\left\{
\begin{array}{ll}
\bbf^j & \text{if } n=0, \\
\bbf^j-\A_1^j \bu_{n-1}^j & \text{otherwise}, 
\end{array}
\right.
\label{eq:matrix-eqn1}
\end{equation}
where $\mathbf{A}_0$ is an $N\times N$ matrix. The set of systems can be recast into a matrix form: 
\begin{equation}
\A_0 \U_n = \B_n,
\label{eq:matrix-eqn}
\end{equation}
where $\U_n=[\bu_n^1, \ldots, \bu_n^J]$ is the solution matrix and $\B_n=[\bb_n^1, \ldots, \bb_n^J]$ contains RHS information from individual problems.

Based on the structure and size of the coefficient matrix $\A_0$, a direct solver or an iterative one can be used to find the solutions. 
Generally speaking, Algorithm~\ref{alg:algorithm2} has several computational advantages:  
(i) $\A_0$ only needs to be accessed once for all the problems; when accessing or generating $\A_0$ represents a major bottleneck of a linear solver, this leads to a significant computational advantage \cite{baker2006improving};  
(ii) A factorization of $\A_0$, if used, only needs to be implemented once in solving all the problems; 
(iii) The multiple RHS vectors in $\B_n$ lead to more efficient matrix-matrix products than matrix-vector products \cite[Section 5.6]{duff2017direct}. Moreover, if a block Krylov subspace algorithm is used in solving \eqref{eq:matrix-eqn}, the multiple RHS vectors would enlarge the search space for minimizing residuals and thus could accelerate the convergence. 
We consider the direct solver in the paper and refer the readers to \cite{ju2019numerical} and reference therein for discussions of block iterative solvers.

A direct solver contains the build phase and the solve phase, in which the former normally has higher complexity than the later. For a dense linear system, the factorization takes $\cO(N^3)$ while the solve phase requires $\cO(N^2)$. Thus, using Algorithm~\ref{alg:algorithm2} is more efficient than solving problems separately as factorization is needed only once.   
Since $\A_0$ is sparse and SPD in the finite element discretization, fast direct solvers are available: using the nested dissection methods \cite{george1973nested} or multi-frontal methods, the complexity for building the LU/Cholesky factorization is $\cO(N^{3/2})$ in 2-D problem and $\cO(N^{2})$ in 3-D (see, for instance,  \cite[Section 9.3]{duff2017direct} and \cite[Section 7.6]{davis2006direct}). Solving the resulting triangular systems has the complexity $\cO(N\log N)$ in 2-D and $\cO(N^{4/3})$ in 3-D. 
We next denote the cost of the build stage by $\cO(N^p)$ and that of the solve stage by $\cO(C_s)$, and analyze the total complexity of Algorithm~\ref{alg:algorithm2} for solving $J$ problems. 
Suppose the algorithm converges after $K-1$ iterations, the total computation complexity is $\cO(N^p+KJC_s)$. For which,  
one factorization of $\A_0$ costs $\cO(N^p)$ and 
$KJ$ times solves of two triangular systems cost $\cO(KJC_s)$. 
But solving $J$ problems separately has the complexity $\cO(J(N^p+C_s))$.   
The ratio of the latter to the former provides a speedup factor  
$$
S_f = \frac{J(N^p+C_s)}{N^p+KJC_s} 
= \frac{1+C_s N^{-p}}{J^{-1}+KC_sN^{-p}}.
$$    
Thus, given a finite $K$, for a fixed $N$, if $J$ is big enough, the factor is on the order of $\frac{N^p}{KC_s}$ that scales as $\cO(\sqrt{N}(\log N)^{-1}K^{-1})$ in 2-D and $\cO(N^{\frac{2}{3}}K^{-1})$ in 3-D; while for a fixed $J$, as $N$ becomes sufficiently large, the factor scales as $J$ in both 2-D and 3-D.

\subsection{Grouping}
Given a group of problems, the convergence of the iterative algorithm could be slow if $\rho$ is close to 1. In such a case, we propose to divide the problems into smaller groups to ensure a small $\rho$ for each group. Motivated by the Centroidal Voronoi tessellations (CVT) method \cite{du1999centroidal}, we introduce the following algorithm. 

Consider a set $W \subset \mathbb{R}^{\widehat{d}}$. A set $\{V_i\}_{i=1}^{n_c}$ is a tessellation of $W$ if $V_i\cap V_j=\emptyset$ for $i\neq j$ and $\cup_{i=1}^{n_c} V_i = W$. Let $|\cdot|$ denote the Euclidean norm on $\mathbb{R}^N$ and, for any $v$ with magnitude greater than $0$, define
$$r(x, v) = \frac{|x-v|}{|v|}.$$
Given a set of points $Z= \{z_i\}_{i=1}^{n_c}$, define the set $V_i$ by 
\begin{equation*}
V_i = \left\{x\in W: r(x, z_i) \leq r(x, z_j) \text{ for } j=1, \ldots, k,\, j\neq i\right\},
\end{equation*} 
where the equality holds only for $i<j$. 
For any region $V\subset \mathbb{R}^{\widehat{d}}$, $z^*$ denotes the center of the region. In a discrete setting, suppose there are $n_s$ elements $\{x_i\}_{i=1}^{n_s}\in V$, we have 
\begin{equation*}
z^* = \frac{1}{n_s} \sum_{i=1}^{n_s} x_i.
\end{equation*} 

Given a set of samples, we need to determine the regions $V_i$ and centers $z_i$. 
Thus, we design an algorithm, presented in Algorithm~\ref{alg:grouping}, that uses iterations to find them. 
First, it initializes the regions (also referred to be clusters) and centers (also referred to be cluster generators). For each data point $x\in W$, the algorithm finds the nearest generator $z\in Z$. 
If that does not match the cluster to which the data point is currently assigned, the 
data point is ``transferred" to the cluster associated with the nearest generator. 
When all data points have been considered, the cluster generators $z$ are
replaced by centers of the clustered data points. 
The process ends when there is no further ``transfers" occurred. 

\begin{algorithm}[!htp]
\small
\SetAlgoLined
\KwIn{ Sample set $W$, number of samples $n_s$, number of centers $n_c$, maximum number of iterations $iter_{max}$.}
\KwOut{ Centers $Z$, Regions $V_1, \ldots, V_{n_c}$.}
Initialize centers $Z=\{z_1, \ldots, z_{n_c}\}$ uniformly spaced in the set\; 
Initialize regions $\{V_1, \ldots, V_{n_c}\}$ associated to each $z_j$ using the empty set, for $j=1, \ldots, n_c$\; 
\For{$k = 1: iter_{\max}$}{
\For{$i=1:n_s$}{
	Take $i$-th sample $x_i = W(i)$ \;
	Determine $z_j = \arg\min\limits_{v\in Z} r(x_i, v)$, where $r(x_i,v) = \frac{|x_i-v|}{|v|}$\;
	Include $x_i$ into the $j$-th region $V_j$ controlled by $z_j$\;
}
Update $z_j$ by the center of $V_j$, for $j=1, \ldots, n_c$\;
\If{no transfers occurred}{Exit the loop}
}
\caption{Grouping algorithm \label{alg:grouping}}
\end{algorithm}

\section{Applications}\label{sec-6}
In this section, we conduct three sets of numerical tests on three application  problems.  The first one solves a parameter-dependent diffusion problem, the second solves a random diffusion problem, and the third simulates a random convection-dominated convection-diffusion problem which is a randomized double-glazing problem \cite{IFISS}. 
Both 1-D and 2-D numerical test results are presented to demonstrate the efficiency 
of the proposed iterative method and to validate the theoretical convergence results.
Particularly, the last two tests show that the proposed iterative method can successfully solve some strongly random diffusion and convection-diffusion problems which the multi-modes method fails to do.

\subsection{Parameter-dependent diffusion problems}\label{sec-6.1}

For a finite set of parameters $\{\omega_j\}_{j=1}^J$, let $a_j(x)=a(\omega_j,x)$, 
$f_j=f(\omega_j,x)$ and consider the following parameter-dependent diffusion problems:
\begin{subequations}\label{eq5.1}
\begin{alignat}{2}
	-\nab\cdot\bigl(a_j\nab u_j\bigr) &= f_j &&\qquad\mbox{in}~ D,\\
	u_j &= 0 &&\qquad \mbox{on}~ \p D,
\end{alignat}
\end{subequations}
where $D$ is a bounded Lipschitz domain in $\mathbb{R}^d$ for $d = 1,2,3$. 

We assume the above PDE is uniformly elliptic,  that is, there exist two constants $0<\lambda< \Lambda$ such that the diffusion coefficient   $\{a_j(x)\}_{j=1}^J$ satisfies  
\begin{align}
\lambda \leq a_j(x) \leq \Lambda \qquad\qquad\forall x\in D, \, 	j = 1,2,\cdots, J.
\end{align}
In addition, assume that $a_j(x) = a_0(x) + \eta_j(x)$ with $\eta_j(x)=\eta(\omega_j,x)$ and $a_0(x)$ satisfies
\begin{enumerate}
	\item[(i)] Uniform ellipticity: $0 <\underline{a}_0 \leq a_0(x) \leq \overline{a}_0$.
	\item[(ii)] Relative dominance: there exists a number $\rho \in (0,1)$ such that
	\begin{align}
		\frac{\|a_j - a_0\|_{L^{\infty}}}{\underline{a}_0} = \frac{\|\eta_j\|_{L^{\infty}}}{\underline{a}_0} < \rho \qquad\forall j = 1,\cdots, J.
	\end{align}
\end{enumerate}

There are a couple of options to choose $a_0(x)$ which are given below (cf. \cite{Yan}). 
\begin{itemize}
	\item[(a)]  $a_0(x)$ is chosen as the arithmetic average of $\{a_j(x)\}_{j=1}^J$, that is,
	\begin{equation}
	a_0(x):= \frac{1}{J} \Bigl(a_1(x)+a_2(x)+\cdots + a_J(x) \Bigr). 
	\end{equation}
	
	\item[(b)]  $a_0(x)$ is chosen as the largest value of $\{a_j(x)\}_{j=1}^J$, that is,
	\begin{equation}
	a_0(x):= \max_{1\leq j\leq J}  a_j(x). 
	\end{equation}
	
\end{itemize}

To fit the abstract framework, we set $V=H^1_0(D)$ and define the bilinear forms $A, A_0, A_1$ and the linear 
functional $F$ respectively as  
\begin{alignat*}{2}
&A(\omega;u,v) = \bigl(a(\omega,\cdot)\nab u, \nab v\bigr), \qquad && F(\omega;v) = \bigl(f(\omega,\cdot),v\bigr),\\
&A_0(u,v) = \bigl(a_0\nab u,\nab v\bigr), \qquad && A_1(\omega;u,v) = \bigl(\eta(\omega,\cdot)\nab u,\nab v\bigr).
\end{alignat*}
It is easy to verify that $A, A_0, A_1$ and $F$ satisfy the convergence criteria laid out in Sections \ref{sec-3} and \ref{sec-4}. Hence, Theorems \ref{thm3.1}, \ref{thm4.2} and \ref{thm4.3} apply to this problem with $\ell=r$. 

Two test cases are considered in this section: one is a 1-D diffusion with an analytic solution, the other is a 2-D diffusion with more degrees of freedom but no analytic solution. We use the former to illustrate the theoretical results and the latter for checking the numerical efficiency. 
To investigate the numerical performance of the iterative algorithm, we consider two metrics: 
$$\mathcal{E}_j = \|u_j-u_{j,n}^h\|_{H^1} \quad \text{  and  }\quad \mathcal{E}_j^h = \|u_j^h-u_{j,n}^h\|_{H^1},$$ 
where $u_j$ is the exact solution, $u_j^h$ is the finite element solution from individual simulations and $u_{j,n}^h$ is the finite element solution at $n$-th iteration of the iterative algorithm.

\medskip
{\bf Test 1.} 
Consider a 1-D diffusion problem (with homogeneous boundary condition) on $D= [0,1]$
with the diffusion coefficient  
$$a_j(x) = 1+x+\epsilon_j \sin(x).$$ 
The exact solution is given by $u(x) = x(x-1)+0.5\sin(20\pi x)+\epsilon_j \sin(40 \pi x)$, and the source term is determined by plugging the exact solution in \eqref{eq5.1}. 
For investigating the iterative algorithm, we consider $J=5$ problems with the choices of $\epsilon_j \in \{0.1035, 0.0727, -0.0303, 0.0294, -0.0787\}$.

For the spatial discretization, the quadratic conforming finite elements ($\ell=2$) are used on a uniform mesh with the size $h$. 
In the iterative algorithm, we set the stopping criterion to be $\max_j \|u_{j,n+1}^h-u_{j,n}^h\|_{H^1}<\mathsf{tol}=10^{-4}$, that is, the maximum $H^1$ norm of the differences between numerical solutions of adjacent iterations is less than $\mathsf{tol}$.

\begin{figure}[htp]
\includegraphics[width=1\linewidth]{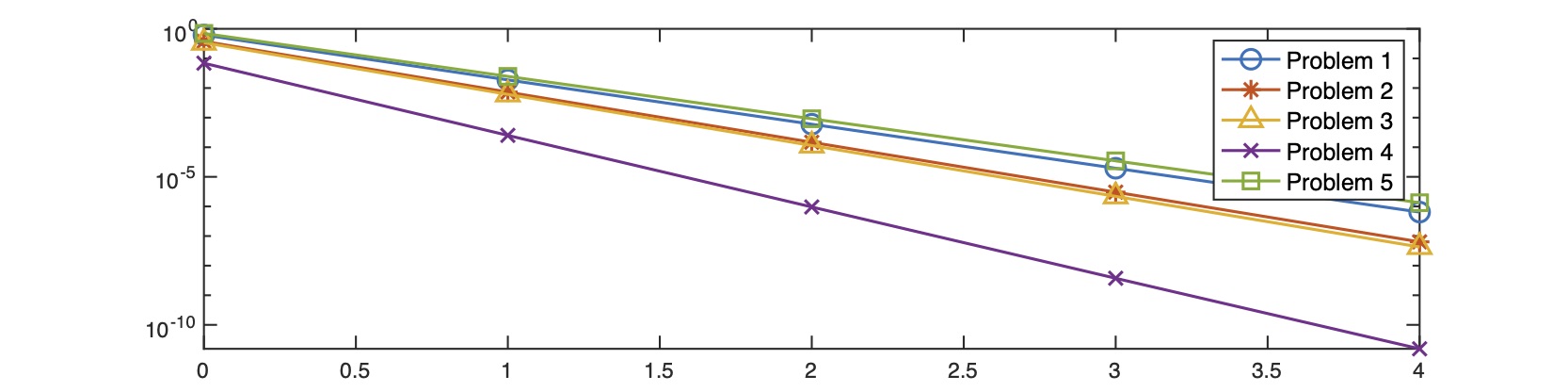}
\caption{Evolution of $\mathcal{E}_j^h$, for $j=1,\ldots, 5$, during the iteration. }
\label{fig:ite1d_1}
\end{figure}

Firstly, we choose $a_0(x) = 1+x+\overline{\epsilon}\sin(x)$ with $\overline{\epsilon}\coloneqq\frac{1}{5}\sum_{j=1}^5 \epsilon_j = 0.0193$. The iterative algorithm takes 4 iterations to complete for meshes with different resolutions. The numerical errors $\mathcal{E}_j$, for $j=1, \ldots, 5$, are listed in Table \ref{tab:err1}, which shows for each of the five problems, the approximation error decays at the second order as mesh is uniformly refined. 
The convergence history of $\mathcal{E}_j^h$ for the five problems is shown in Figure \ref{fig:ite1d_1} by taking the $h=2^{-10}$ case for example. 
Applying regressions on the results shows 
$\mathcal{E}_1^h\sim \mathcal{O}(0.032^{n+1})$ in problem 1, 
$\mathcal{E}_2^h\sim \mathcal{O}(0.021^{n+1})$ in problem 2, 
$\mathcal{E}_3^h\sim \mathcal{O}(0.019^{n+1})$ in problem 3, 
$\mathcal{E}_4^h\sim \mathcal{O}(0.0039^{n+1})$ in problem 4, and 
$\mathcal{E}_5^h\sim \mathcal{O}(0.038^{n+1})$ in problem 5. 
The ratios $\rho_j = \frac{\|a_j-a_0\|_{L^\infty}}{\underline{a}_0}$ for these five problems are $0.071$, $0.045$, $0.042$, $0.009$, and $0.083$, respectively, which are close and proportional to the regression rates. 
We found a better ratio that matches the regression analysis result can be defined by $\widehat{\rho}_j = \left\|\frac{a_j-a_0}{a_0}\right\|_{L^\infty}$, whose values are respectively 
$0.035$, $0.022$, $0.021$, $0.004$ and $0.041$ for these problems. Further analysis in this aspect will be performed in a future work.

\begin{table}[htp]\small
\renewcommand{\arraystretch}{1.4}
\centering
\begin{tabular}{ |c|c|c|c|c|c| } 
 \hline
           $h$      & $\mathcal{E}_1$ & $\mathcal{E}_2$  &  $\mathcal{E}_3$  &  $\mathcal{E}_4$  &  $\mathcal{E}_5$   \\ 
 \hline
 $\frac{1}{2^7}$    & $3.82\times 10^{-1}$ & $3.03\times 10^{-1}$ & $2.21\times 10^{-1}$ & $2.19\times 10^{-1}$ & $3.18\times 10^{-1}$ \\ 
 $\frac{1}{2^8}$    & $9.62\times 10^{-2}$ & $7.63\times 10^{-2}$ & $5.54\times 10^{-2}$ & $5.50\times 10^{-2}$ & $8.00\times 10^{-2}$  \\ 
 $\frac{1}{2^9}$    & $2.41\times 10^{-2}$ & $1.91\times 10^{-2}$ & $1.39\times 10^{-2}$ & $1.38\times 10^{-2}$ & $2.00\times 10^{-2}$  \\ 
 $\frac{1}{2^{10}}$ & $6.03\times 10^{-3}$ & $4.78\times 10^{-3}$ & $3.46\times 10^{-3}$ & $3.44\times 10^{-3}$ & $5.01\times 10^{-3}$  \\ 
 \hline
 \end{tabular}
 \caption{Errors $\mathcal{E}_j$, for $j=1, \ldots, 5$, of the iterative algorithm solutions at different $h$.}
 \label{tab:err1}
\end{table}

\begin{figure}[htp]
\includegraphics[width=1\linewidth]{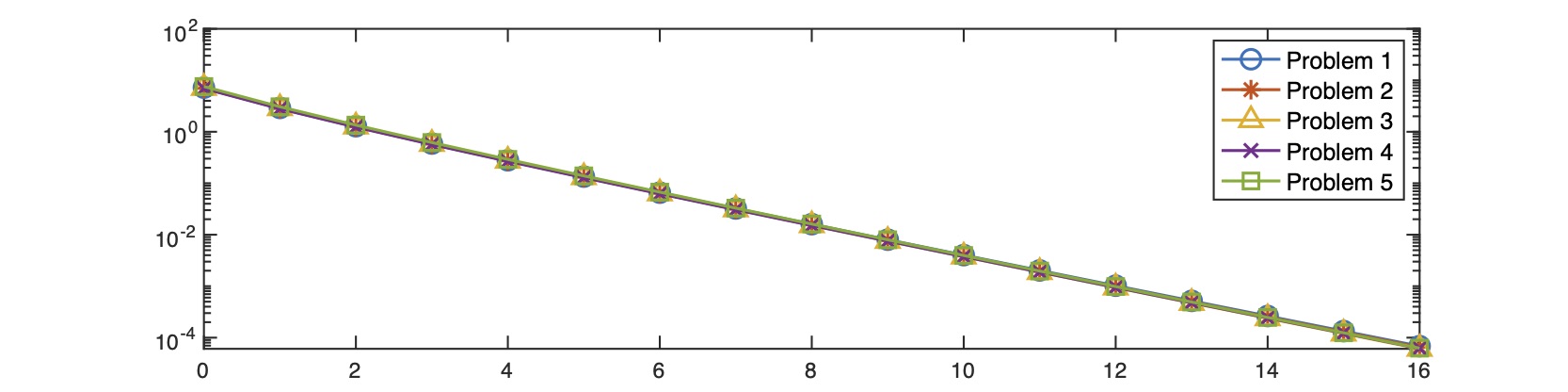}
\caption{Evolution of $\mathcal{E}_j^h$, for $j=1, \ldots, 5$, during the iteration. }
\label{fig:ite1d_2}
\end{figure}

Secondly, we choose $a_0 = a_{\infty} \coloneqq \max_j \max_{x\in D} |a_j(x)| = 2.0871$ while keeping the same computational setting as the previous case. 
The iterative algorithm achieves the same errors listed in Table \ref{tab:err1}, which is not surprising because the same stopping criterion is used in both tests. However, the total number of iterations increases to 16 in this case. 
The evolution of $\mathcal{E}_j^h$ for the five problems during the iteration when $h=2^{-10}$ is plotted in Figure \ref{fig:ite1d_2}. 
Regressions on the data shows $\mathcal{E}_j^h\sim \mathcal{O}(0.49^{n+1})$ in these 5 problems. 
Note that the ratios $\rho_j=\widehat{\rho}_j$ in this case because $a_0$ is a constant. For the five problems, the ratios are about 
$0.52$ that are close to the regression results. 

It is observed that the numerical results match the theoretical analysis in Theorems~\ref{thm4.2} - \ref{thm4.3}. 
In particular, $\mathcal{E}_j$ is dominated by the finite element discretization error, that is of the order $\mathcal{O}(h^\ell)$ with $\ell=2$; and $\mathcal{E}_j^h$ is of the order $\mathcal{O}(\rho_j^{n+1})$, for $j=1, \ldots, 5$. 

\medskip
{\bf Test 2.} 
We consider a diffusion problem on the domain $D=[-1,1]^2$. 
Let $D_1$ be a disk centered at the origin of radius 0.5 and $D_0 = D\backslash \overline{D}_1$. 
The conductivity $a(x)$ is a piecewise constant on $D$: 
$$
a(x)|_{D_1} = \mu^{[1]} \quad \text{and}\quad a(x)|_{D_0} = 1.
$$
The problem has zero forcing and is associated to a Dirichlet boundary condition on the top edge and Neumann boundary conditions on the other sides. In particular, 
\begin{align*}
u = 0 & \text{ on } \Gamma_{top}= [-1, 1]\times \{1\}, \\
a\nabla u \cdot n = 0 &\text{ on } \Gamma_{side} = \{\pm 1\}\times (-1,1), \\
a\nabla u \cdot n = \mu^{[2]} &\text{ on } \Gamma_{bottom} = [-1, 1]\times \{-1\}.
\end{align*}
Denote the parameter vector by $\boldsymbol{\mu} = (\mu^{[1]}, \mu^{[2]})$, whose ranges is specified by 
\[
\Ome:=\{\boldsymbol{\mu} \in \mathbb{R}^2; \, \, 0.1\leq \mu^{[1]} \leq 10, \, -1\leq \mu^{[2]} \leq 1 \}.
\]
\begin{figure}[htp]
\centering
\begin{minipage}{0.46\textwidth}
	\includegraphics[width=1\linewidth]{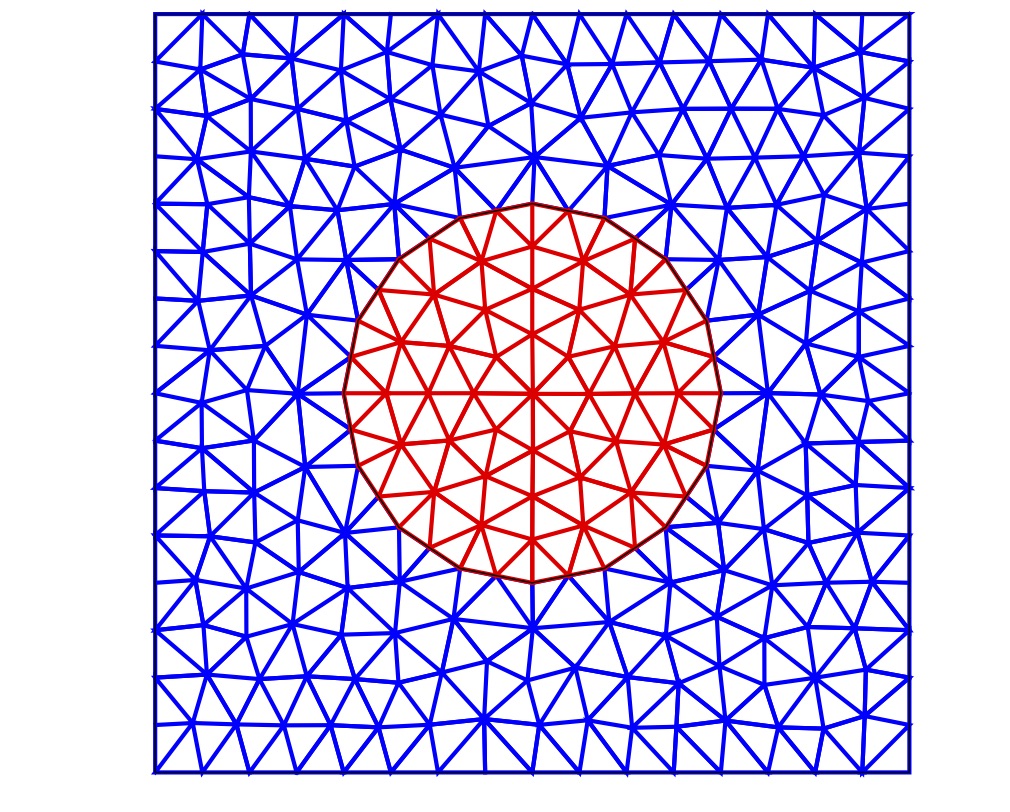}
\end{minipage}
\begin{minipage}{0.46\textwidth}
	\includegraphics[width=1\linewidth]{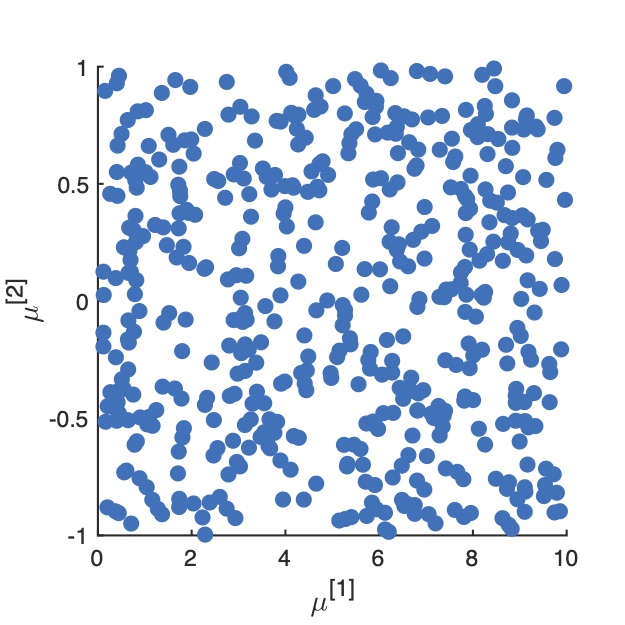}
\end{minipage}
\caption{(Left) A finite element mesh  with 496 elements; (Right) $\boldsymbol{\mu}$-samples.}
\label{fig:para1}
\end{figure}
In this case, the bilinear forms $A, A_0, A_1$ and the linear functional $F$ have the following form: 
\begin{alignat*}{2}
&A(\boldsymbol{\mu};u,v) = \bigl(\mu^{[1]}\nab u, \nab v\bigr)_{D_1}+\bigl(\nab u, \nab v\bigr)_{D_0}, \quad && F(\boldsymbol{\mu};v) = \bigl(f(\boldsymbol{\mu}),v\bigr)+\langle\mu^{[2]}, v\rangle_{\Gamma_{\mbox{\tiny bottom}}},\\
&A_0(u,v) = \bigl(\mu_0^{[1]}\nab u, \nab v\bigr)_{D_1}+\bigl(\nab u, \nab v\bigr)_{D_0}, \quad && A_1(\boldsymbol{\mu};u,v) = \bigl((\mu^{[1]}-\mu_0^{[1]})\nab u, \nab v\bigr)_{D_1}, 
\end{alignat*}
where we use $a_0|_{D_1}\coloneqq \mu_0^{[1]}$ and $a_0|_{D_0}\coloneqq 1$. Since these forms have affine parameter dependence, the assembly of the associated matrices and vectors can be greatly simplified in computation. 

We use the quadratic conforming finite element method for spatial discretization on a triangulation of $D$ with the mesh size $h$. For instance, the triangulation of $h=1/8$ is shown in Figure \ref{fig:para1} (left), which contains $496$ elements in total. In the iterative algorithm, the stopping criterion $\mathsf{tol} = 10^{-4}$ is selected. 

There are $n_s$ $\boldsymbol{\mu}$-samples $\{\boldsymbol{\mu}_j= ( \mu_j^{[1]}, \mu_j^{[2]}) \}_{j=1}^{n_s}$ randomly chosen in the parameter domain, e.g., $n_s= 500$ randomly selected samples are shown in Figure \ref{fig:para1} (right). 
Note that only $\mu_j^{[1]}$ appears in the bilinear form, the ratio $\rho_j$ only depends on $\mu_j^{[1]}$. 
In such a case, choosing $\mu_0^{[1]} \coloneqq \max_j \mu_j^{[1]}$ requires many iterations as shown in previous test, but choosing $\mu_0^{[1]}\coloneqq \frac{1}{J} \sum_{j=1}^{n_s} \mu_j^{[1]}$ also will lead to a big $\rho_j$ for some problems and require many iterations as well. Therefore, we first divide the sample set 
into multiple subgroups using Algorithm~\ref{alg:grouping} and then apply the iterative algorithm to each of them. 

\begin{figure}[htp]
\centering
\begin{minipage}{0.48\textwidth}
	\includegraphics[width=1\linewidth]{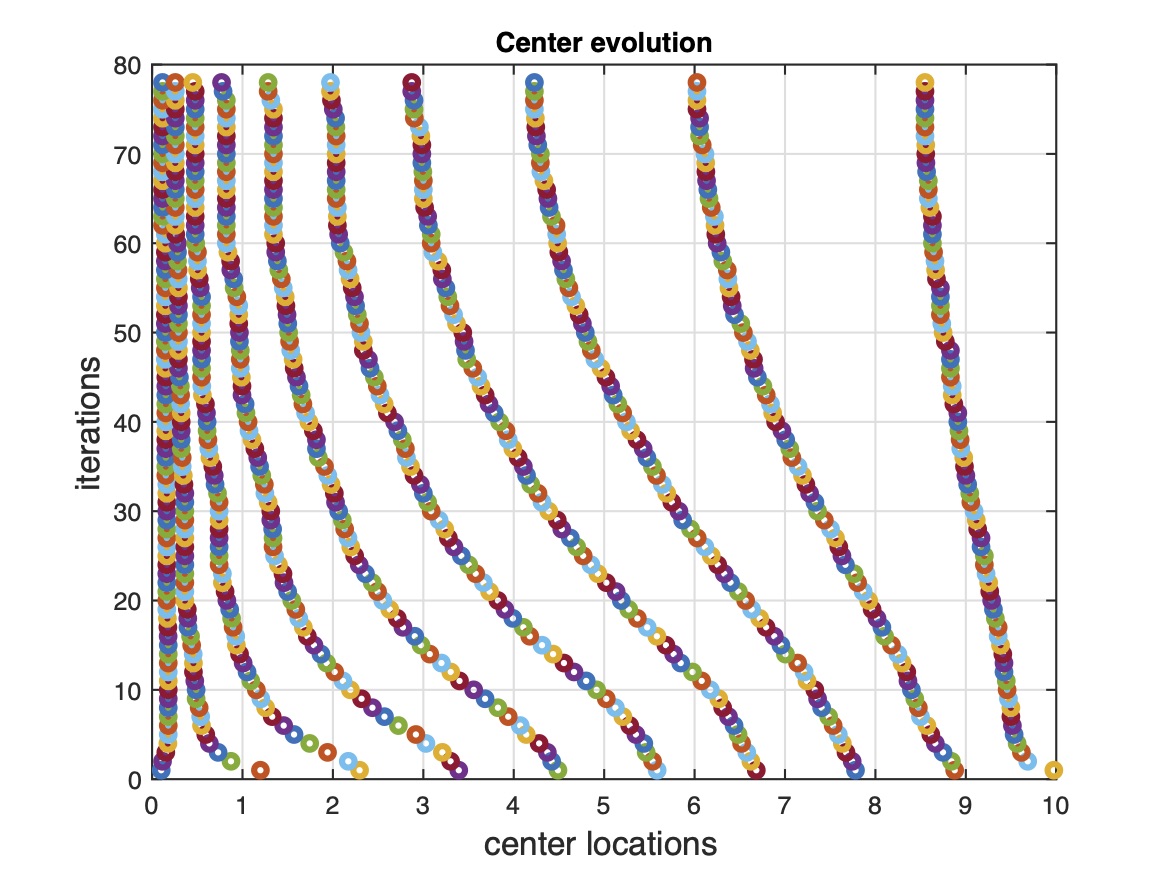}
\end{minipage}
\begin{minipage}{0.48\textwidth}
	\includegraphics[width=1\linewidth]{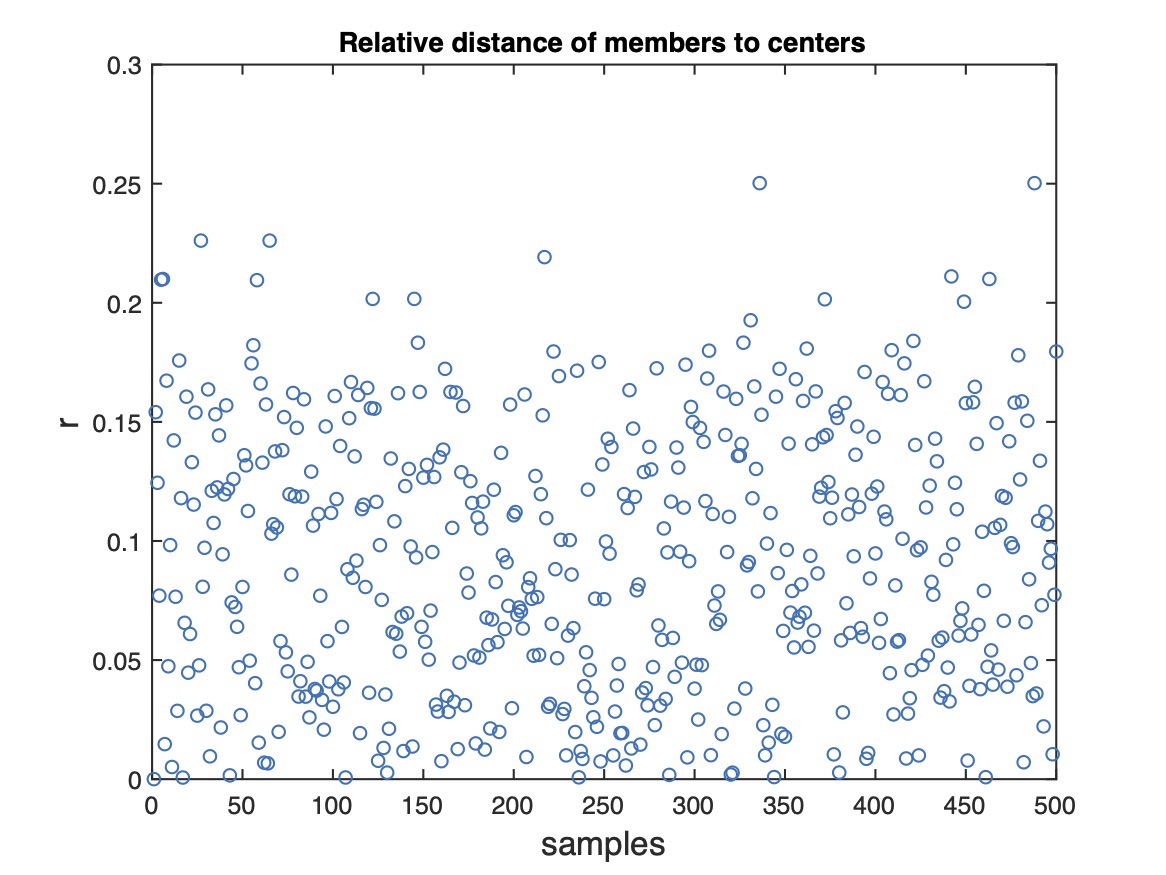}
\end{minipage}
\caption{(Left) Evolution of centers of 10 subgroups; (Right) $r(\mu_i^{[1]}, z_j)$ for each $\boldsymbol{\mu}_i\in W$ belongs to $V_j$ and $z_j$ is the center of $V_j$.}
\label{fig:group}
\end{figure}

\begin{figure}[htp]
\centering
	\includegraphics[width=1\linewidth]{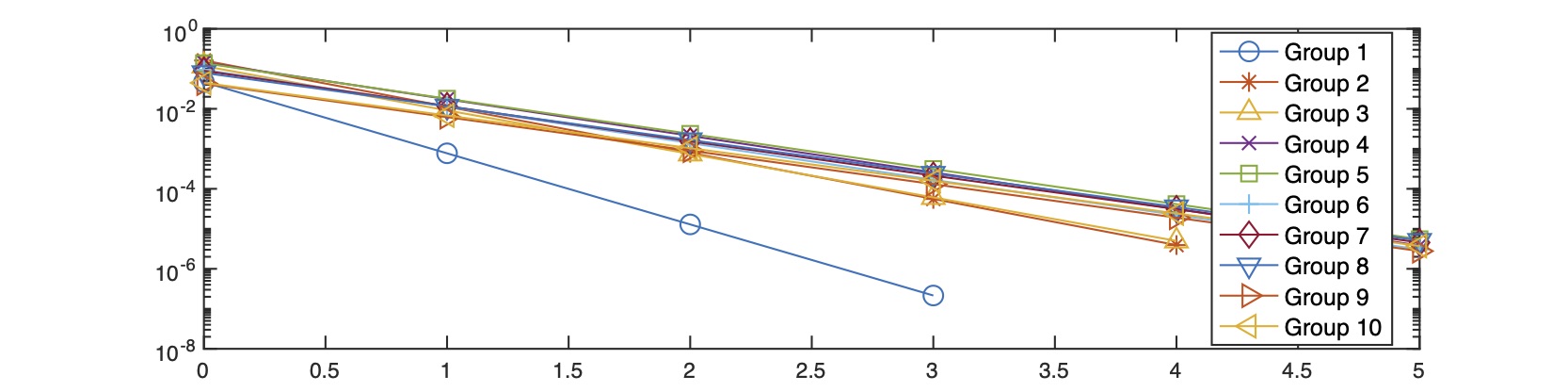}
\caption{Evolution of $\max_j\mathcal{E}_j^h$ for each group during the iteration.}
\label{fig:group2}
\end{figure}

We use a mesh of size $h=\frac{1}{32}$ for spatial discretization, which contains 8,136 elements and 16,529 nodes. 
For the same set including $n_s= 500$ samples, Algorithm~\ref{alg:grouping} is fed with $W\coloneqq \{\mu_1^{[1]}, \ldots, \mu_{n_s}^{[1]}\}$ and finds $n_c=10$ regions and associated centers. 
The evolution of centers during the process is shown in Figure~\ref{fig:group} (left), which shows the sample points stop transferring after 78 steps.  
Each sample from the set is assigned to a region $V_k$ with the center $z_k$, for $k=1, \ldots n_c$. The ratio $r({\mu}_i^{[1]}, z_k)$ for each $\boldsymbol{\mu}_i$ is shown in Figure~\ref{fig:group} (right). It is seen that all these ratios are below 0.3.  
The iterative algorithm is sequentially performed in the $n_c$ groups with $a_0|_{D_1}\coloneqq z_k$ for each group. The entire simulation is completed in 29.48 seconds and the maximum number of iterations for each group is 5.
The detailed information for all the groups is listed in Table \ref{tab:group1}, where $\rho$ represents the maximum value of $r(\mu_i^{[1]}, z_j)$ for elements in $j$-th group. 
The maximum errors between the iterative solution and the finite element solution, $\max_j \mathcal{E}_j^h$, in each group are plotted in Figure \ref{fig:group2}. 

\begin{table}[htp]\small
\renewcommand{\arraystretch}{1.4}
\centering
\begin{tabular}{ |c|c|c|c|c|c|c|c| } 
\hline
group & size  & $\mu^{[1]}$-region         & center & $\rho$  & iterations & max $\mathcal{E}_j^h$ & time  \\
\hline
1           &    4         & [0.10, 0.13] & 0.12  & 0.12  & 3         & $2.15\times 10^{-7}$     & 0.49 \\
2           &   15         & [0.19, 0.32] & 0.26  & 0.25  & 4         & $3.92\times 10^{-6}$     & 1.06 \\
3           &   14         & [0.36, 0.54] & 0.45  & 0.20  & 4         & $4.99\times 10^{-6}$     & 1.05 \\
4           &   17         & [0.59, 0.94] & 0.77  & 0.23  & 5         & $3.81\times 10^{-6}$     & 1.37 \\
5           &   33         & [1.01, 1.55] & 1.28  & 0.21  & 5         & $5.48\times 10^{-6}$     & 2.18 \\
6           &   32         & [1.61, 2.33] & 1.97  & 0.18  & 5         & $3.05\times 10^{-6}$     & 2.12 \\
7           &   50         & [2.36, 3.39] & 2.87  & 0.18  & 5         & $4.55\times 10^{-6}$     & 3.07 \\
8           &   68         & [3.49, 4.97] & 4.23  & 0.17  & 5         & $5.08\times 10^{-6}$     & 3.90 \\
9           &  111         & [4.99, 7.06] & 6.03  & 0.17  & 5         & $2.78\times 10^{-6}$     & 6.15 \\
10          &  156         & [7.12, 9.97] & 8.55  & 0.17  & 5         & $3.72\times 10^{-6}$     & 8.09 \\
 \hline
 \end{tabular}
 \caption{Information of subgroups: size, region, center and maximum $\rho$; performance of the iterative algorithm including number of iterations, maximum errors and wall-clock time for simulations.}
 \label{tab:group1}
\end{table}

Next, we vary the number of samples and compare the iterative algorithm with the individual simulations in terms of wall-clock times. Since the exact solution is unknown, we measure the maximum difference in $H^1$ norm between iterative algorithm solutions and individual solutions. The results are listed in Table \ref{tab:group2}, where $T_{it}$ represents the time elapsed for integration in the iterative algorithm and $T_{ind}$ is that for integrating individual systems. It is observed that as the size of sample set increases, the efficiency of iterative algorithm improves. When $n_s=2500$, the iterative algorithm saves about 50\% simulation time. 
\begin{table}[!htp]\small
\renewcommand{\arraystretch}{1.4}
\centering
\begin{tabular}{ |l|l|l|l| } 
\hline
$n_s$   & $T_{it}$  &   $T_{ind}$  &  $\max_j \mathcal{E}_j^h$ \\
\hline
100     &    8.09    &    8.74      &      $ 7.2\times 10^{-6}$         \\
500     &   29.48    &    43.14     &      $ 5.48\times 10^{-6}$         \\
2500    &   116.10   &   233.23     &      $ 1.27\times 10^{-6}$         \\
\hline
\end{tabular}
\caption{Comparison of the iterative algorithm with individual simulations at $h=1/32$: wall-clock times and maximum values of $\mathcal{E}_j^h$ for $j=1, \ldots n_s$.}
\label{tab:group2}
\end{table}

Finally, we fix $n_s=2500$ while doubling the number of centers $n_c$ from 5 to 160. The wall-clock times $T_{it}$ are listed in Table \ref{tab:group3}, which shows the simulation time first decreases then increases. The least computational time is achieved when $n_c=80$, which saves over 60\% simulation time compared to the individual simulations. It is also observed that as the number of groups increases, the maximum difference between iterative algorithm and individual simulation solutions, $\max_j \mathcal{E}_j^h$, decreases, which is because $\rho$ shrinks as the problems are divided into more groups.
 
\begin{table}[!htp]\small
\renewcommand{\arraystretch}{1.4}
\centering
\begin{tabular}{ |l|l|l|l|l|l|l| } 
\hline
$n_c$     & 5       &  10      &  20      &  40     &  80      &  160      \\
\hline 
$T_{it}$  & 187.5  &  116.1  &  102.64  &  93.01  &  89.69   &  111.75   \\
$\max_j \mathcal{E}_j^h$ &  $0.0000551$ &  $0.0000127$ &  $0.000005$ &  $0.00000327$ &  $0.00000328$ &   $0.00000167$ \\   
\hline
\end{tabular}
\caption{Wall-clock times and maximum values of $\mathcal{E}_j^h$ in the iterative algorithm for $2500$ samples that are divided into $n_c$ groups.}
\label{tab:group3}
\end{table}

\subsection{Random diffusion equations}\label{sec-6.2}
In this subsection we consider the following random diffusion test problem:  
\begin{subequations}\label{eq5.5}
\begin{alignat}{2}
	-\nab\cdot\bigl(a(\omega,\cdot)\nab u\bigr) &= f(\omega,\cdot) &&\qquad \mbox{in}~ D, \quad\mbox{$\mP$-a.s. }\\
	u &= 0  &&\qquad \mbox{on}~ \p D, \quad  \,\mbox{$\mP$-a.s. } 
\end{alignat}
\end{subequations}
Where $D$ is a bounded Lipschitz domain in $\mathbb{R}^d$ for $d = 1,2,3$. $a(\omega,x)$ and $f(\omega,x)$ are random fields on the probability space $(\Omega, \mathbb{P}, \mathcal{F})$ and they 
are assume to have continuous and bounded covariance function. $\mathbb{E}[\cdot]$ denotes 
the expectation operator.  In addition, $a$ satisfies $\mP$-a.s.
the following condition:
\begin{align}
	0<\lambda \leq a(\omega,x) \leq \Lambda\qquad\forall x\in D.
\end{align}
Moreover, we assume that $\mathbb{E}[\|f\|_{L^2(D)}]<\infty$.
Below we propose two approaches for utilizing the abstract framework for problem \eqref{eq5.5}. 
The first approach mimics the multi-modes method of \cite{Feng} while the second one mimics
the ensemble method of \cite{Yan}. 

\subsubsection{Approach \#1}
Suppose that there exists following decomposition of $a$:
\begin{align}
a(\omega,x) = a_0(x) + \eta(\omega,x) \qquad\mbox{a.s.}~\forall x \in D,
\end{align}
where $a_0$ is independent of $\omega$ and also satisfies 
\begin{enumerate}
	\item[(i)] {\em Uniform ellipticity:} $0 <\underline{a}_0 \leq a_0(x) \leq \overline{a}_0$ for all $x \in D$. 
	\item[(ii)] {\em Relative dominance:} there exists a number $\rho \in (0,1)$ such that
	\begin{align}
	\mP\bigl\{\omega \in \Ome;\, \|\eta(\omega,\cdot)\|_{L^{\infty}(D)} \leq \rho\, \underline{a}_0 \bigr\} = 1.
	\end{align}
\end{enumerate}

To apply the abstract framework, we define 
\begin{alignat*}{2}
&A(\omega; u,v) = \bigl(a(\omega,\cdot)\nab u, \nab v\bigr), &&\qquad F(\omega;v) = \bigl(f(\omega,\cdot),v\bigr),\\
&A_0(u,v) = \bigl(a_0\nab u,\nab v\bigr), &&\qquad 
A_1(\omega; u,v) = \bigl(\eta(\omega,\cdot)\nab u,\nab v\bigr).
\end{alignat*}
It is easy to verify that $A, A_0, A_1$ and $F$ satisfy $\mP$-a.s. the convergent criteria laid out in Sections \ref{sec-3} and \ref{sec-4}. Hence, Theorems \ref{thm3.1}, \ref{thm4.2} and \ref{thm4.3} apply to this problem with $\ell=r$.  In particular, there holds
\begin{align}\label{rate_approach2}
\mE\bigl[\|u - U^h_n\|_{H^1}\bigr] \leq C\bigl(h^r + \rho ^{n+1}\bigr).
\end{align}
It should be noted that $C>0$ is independent of $\omega$  because 
$\mathbb{E}[\|f\|_{L^2(D)}]<\infty$. 

We remark that in this first approach we do not specify how to discretize the stochastic variable $\omega$. However, since the abstract framework suits best with the random sampling method, in practice,  the Monte Carlo method would be used to discretize the right-hand sides in Step 1 and 2 of Algorithm 2.  

Finally, we note that the choice of $a_0$ depends on the structure of $a$ (cf. \cite{Feng}). 

To verify the convergence rates in \eqref{rate_approach2} of this approach, we consider the following random coefficient boundary value problem \cite{Feng} for our numerical tests:
\begin{align}\label{boundary_problem2}
-\frac{d}{dx}\Bigl(\bigl(1+\varepsilon X(\omega)\bigr)\frac{d u(\omega,x)}{dx}\Bigr) &= X(\omega), \qquad 0 < x < 1,\\\nonumber
u(\omega,0) =0,\quad u(\omega,1) &= 0,
\end{align}
where $X(\omega)$ is a uniformly distributed random variable defined in 
a probability space $(\Omega, \mathcal{B}, \mathbb{P})$, where the sample space $\Omega=[0,1]$,
$\mathcal{B}$ denotes the $\sigma$-algebra of the Borel sets and $\mathbb{P}$ denotes the 
Lebesgue probability measure. In addition, we assume that $\varepsilon>0$. 
The true solution of \eqref{boundary_problem2} is given by 
\begin{align}\label{exact_sol}
u(\omega,x) = \frac{X(\omega)}{2(1+\varepsilon X(\omega))}(x-x^2).
\end{align}

We have that $a(\omega,x) = 1+\varepsilon X(\omega)$ satisfies the elliptic condition with $\lambda = 1, \Lambda = 1+\varepsilon$. In addition, $a(\omega,x) = a_0(x) + \eta(\omega,x)$, where $a_0(x) = 1$ and $\eta(\omega,x) = \varepsilon X(\omega)$ satisfies the conditions (i) with $\underline{a}_0 = 1$ and (ii) with $\rho = \varepsilon$. In order to fit the abstract framework, in this example, we define $V = H^1_0(0,1)$ and
\begin{alignat*}{2}
&A(\omega; u,v) = \bigl((1+\varepsilon X(\omega)) u',  v'\bigr), &&\qquad F(\omega;v) = \bigl(X(\omega),v\bigr),\\
&A_0(u,v) = \bigl( u', v'\bigr), &&\qquad 
A_1(\omega; u,v) = \bigl(\varepsilon X(\omega) u', v'\bigr).
\end{alignat*}

Let $V^h_1$ denote the standard linear finite element subspace of $V$ which is used for the spatial discretization of \eqref{boundary_problem2}. Let $\{U^h_n\}_{n\geq 0}$ be the approximate solution from Algorithm 2, where $n$ is the number of iterations. Let 
\begin{align}
\mathcal{E}^n_{L^2} := \mE\bigl[\|u - U^h_n\|_{L^2}\bigr],\quad \mathcal{E}^n_{H^1} := \mE\bigl[\|u - U^h_n\|_{H^1}\bigr].
\end{align}

As mentioned above, we use the Monte Carlo method to discretize the stochastic variable 
$\omega$ on the right-hand side in Step 1 and 2 of Algorithm 2. The number of the Monte Carlo samples is chosen to be $J = 10^4$. To test the validity and accuracy of the proposed iterative method, we set $h = 0.01$ for spatial discretization and then vary the number of iterations $n$.
Table \ref{tab6.7} displays the $H^1$-norm errors of the iterative solutions for different values of $n$ and $\varepsilon$. As expected,  the computed solutions have smaller errors for
smaller $\rho=\varepsilon$, and the method converges as long as $\varepsilon<1$ which is predicted
by our convergence theorem. In addition, our numerical tests show that the method ceases to 
converge for $\rho=\varepsilon>1$, this indicates that our convergence result is sharp.

\begin{table}[htp]\small
	\renewcommand{\arraystretch}{1.4}
	\centering
	\begin{tabular}{ |l|l|l|l|l|l|l| } 
		\hline
		$\rho$      & $n = 1$ & $n=2$  & $n =3$  & $n =4$ &$n=5$ &$n=6$ \\ 
		\hline 
		$0.4$    & $0.093$ &$0.0033$  &$0.0016$  & $0.00123$  &$0.00115$ &$0.00114$\\ 
		$0.6$    & $0.0184$ &$0.0089$  & $0.0047$  &$0.0027$  &$0.0018$ & $0.0014$\\ 
		$0.8$ & $0.0295$ & $0.0187$ & $0.0124$ & $0.0086$  & $0.0062$ &$0.0046$\\ 
		$0.9$ &$0.0356$ & $0.0253$ & $0.188$ &$0.0145$ & $0.0115$ &$0.0093 $\\
		\hline
	\end{tabular}
	\caption{Approach \#1: Errors $\mathcal{E}_{H^1}^n$ of the computed solutions with different $\rho$ and $n$.}
	\label{tab6.7}
\end{table}

 To test the convergence orders for the spatial discretization, we fix $\rho =\varepsilon = 0.1$ and the number of iterations $n = 10$, then vary the mesh size $h$. Table \ref{tab6.5} shows 
 the $H^1$- and $L^2$-errors of the computed solution by Algorithm 2.
\begin{table}[htp]\small
	\renewcommand{\arraystretch}{1.4}
	\centering
	\begin{tabular}{ |l|l|l|l|l| } 
		\hline
		$h$      & $\mathcal{E}_{H^1}^n$ & Order  & $\mathcal{E}_{L^2}^{n}$  & Order   \\ 
		\hline
		$0.2$    &$2.43\times 10^{-2}$  &  & $1.5000 \times 10^{-3}$ &   \\ 
		$0.1$    & $1.29\times 10^{-2}$ &$0.92$  &$4.0659\times 10^{-4}$  & $1.88$  \\ 
		$0.05$    &$6.60\times 10^{-3}$  &$0.97$  & $1.0443\times 10^{-4}$  &$1.96$  \\ 
		$0.025$ & $3.30\times 10^{-3}$ &$1.00$  &$2.6450\times 10^{-5}$  &$1.98$   \\ 
		\hline
	\end{tabular}
	\caption{Errors $\mathcal{E}_{H^1}^n$ and $\mathcal{E}_{L^2}^n$, for $n =10$, of the iterative algorithm solutions at different $h$.}
	\label{tab6.5}
\end{table}
We observe that an $O(h)$ order convergent rate in the $H^1$-norm and an $O(h^2)$ order rate in the $L^2$-norm are obtained for the numerical solution $U^h_n$ which is consistent with the theoretical error estimate in \eqref{rate_approach2}. 

\begin{figure}[htp]
	\begin{center}
		\includegraphics[scale=0.26]{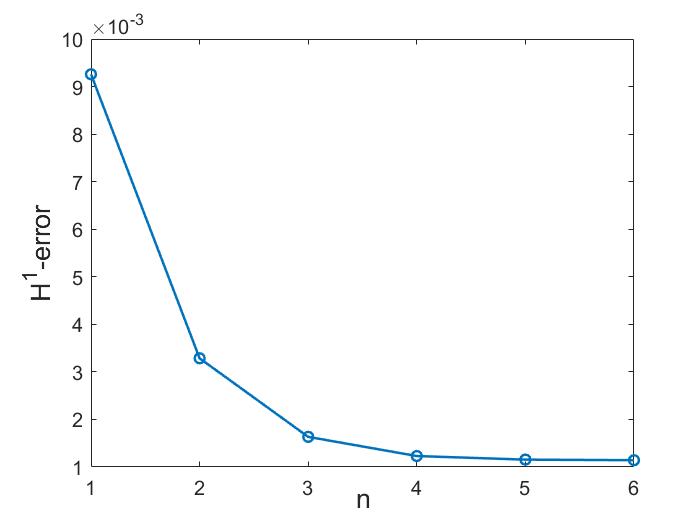}
		\includegraphics[scale=0.26]{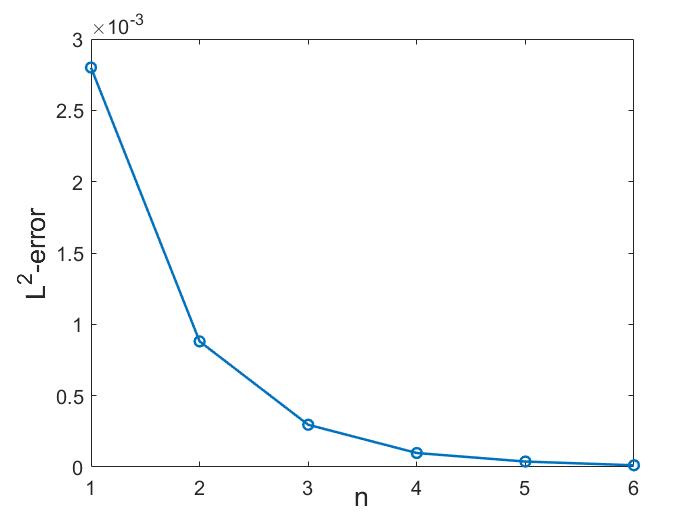}
		\caption{Approach \#1: Errors $\mathcal{E}_{H^1}^n$ (left) and $\mathcal{E}_{L^2}^n$ (right) for $\rho = \varepsilon =0.4,\, h =0.01$, {\tt tol} $= 10^{-4}$ and different $n$.}\label{H1_L2_err}
	\end{center}
\end{figure}

Next, we verify the dependence of the errors on the parameter $\rho=\varepsilon$. 
To the end,  we fix $\rho = 0.4, h = 0.01$ and then vary the number of iterations $n$. We also 
use the stopping criteria $\mE\bigl[\|U^h_{n+1} - U^h_n\|_{H^1}\bigr] < {\tt tol} := 10^{-4}$. 
The $H^1$- and $L^2$-norm errors of the computed solutions are shown in Figure \ref{H1_L2_err}. 
It is clear that the errors are decreasing as the number of iterations increases and 
the convergence order $O(\rho^{n+1})$ is indeed observed as predicted in \eqref{rate_approach2}.
In addition, Figure \ref{H1_L2_err} also shows that only $6$ iterations is needed to 
trigger the stopping criterion.

\begin{figure}[htp]
	\begin{center}
		\includegraphics[scale=0.26]{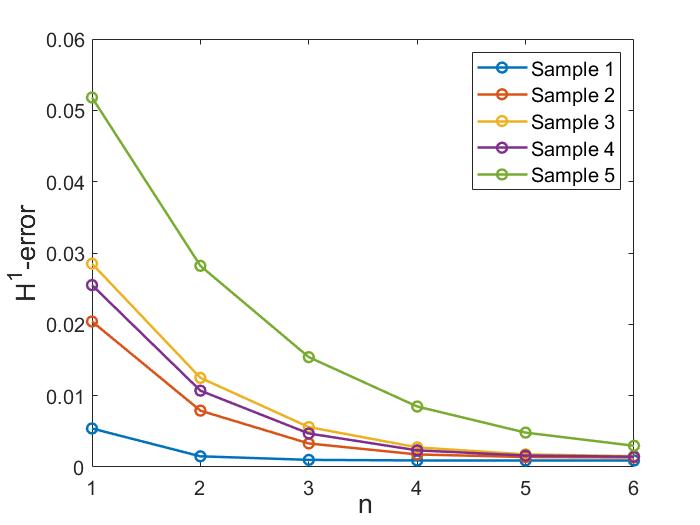}
		\includegraphics[scale=0.26]{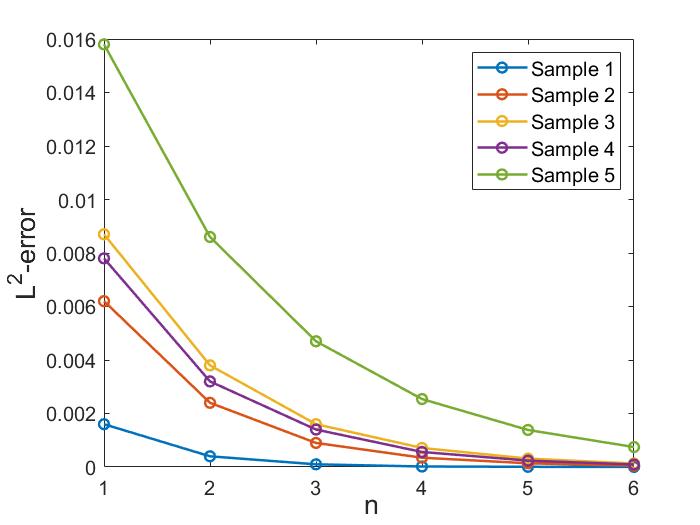}
		\caption{Approach \#1: Errors $\mathcal{E}_{H^1}^n$ (left) and $\mathcal{E}_{L^2}^n$ (right) decay for the first $5$ Monte Carlo samples with $\rho = \varepsilon =0.4,\, h =0.01$ and different $n$.}\label{5samples_1app}
	\end{center}
\end{figure}

Furthermore, we plot the errors of the first five Monte Carlo samples in Figure \ref{5samples_1app} to show that the errors for different samples are different and they are also different from the expected value which is presented in Figure \ref{H1_L2_err}, but they all are decreasing as $n$ increases provided that $\varepsilon$ is small.

\begin{figure}[htp]
	\begin{center}
		\includegraphics[scale=0.26]{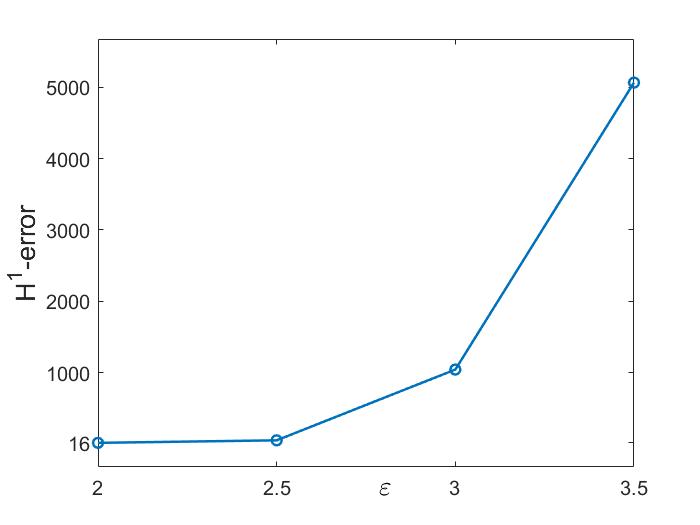}
		\includegraphics[scale=0.26]{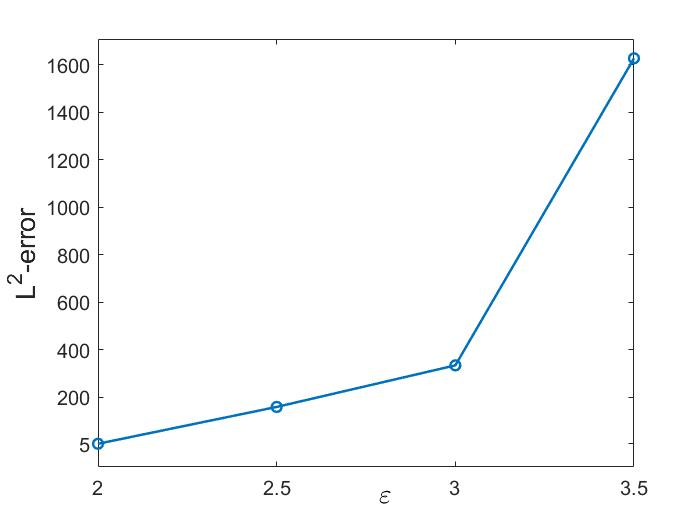}
		\caption{Approach \#1: Errors $\mathcal{E}_{H^1}^n$ (left) and $\mathcal{E}_{L^2}^n$ (right) for different $\varepsilon > 1$ and $n=10$.}\label{relative_dominant}
	\end{center}
\end{figure}

Finally, to check the dependence of the iterative method on the parameter $\varepsilon$, we fix $n = 10,\, h = 0.01$ and then consider $\varepsilon = 2, 2.5, 3, 3.5$, which all are larger than $1$. We note that in these cases, the relative dominant condition (ii) is violated, hence, our convergence results do not apply anymore. 
Figure \ref{relative_dominant} shows evidently that the $H^1$- and $L^2$-norm errors 
increase significantly as $\varepsilon$ becomes larger, which clearly indicates that 
the iterative method may diverge when $\varepsilon>1$ and also shows that the 
relative dominant condition (ii) is sharp. 

\subsubsection{Approach \#2}
In this subsection, we present a different approach to utilize the proposed iterative method for
solving the same random diffusion problem \eqref{eq5.5} based on the Monte Carlo finite 
element discretization. To solve \eqref{eq5.5} by the Monte Carlo method, let $\{\omega_j\}_{j=1}^J$ be $J$ samples 
from the sample space $\Ome$, we then consider (deterministic) problems for $j=1,2,\cdots, J$
\begin{subequations}\label{eq5.6}
	\begin{alignat}{2}
	-\nab\cdot\bigl(a(\omega_j,\cdot)\nab u\bigr) &= f(\omega_j,\cdot) &&\qquad \mbox{in}~ D,\\
	u &= 0 &&\qquad  \mbox{on}~ \p D. 
	\end{alignat}
\end{subequations}
Hence, the Monte Carlo method converts the random diffusion problem \eqref{eq5.5} into the parameter-dependent 
problem \eqref{eq5.6} which is of the same type as problem \eqref{eq5.1}. As a result, the method described in Section \ref{sec-6.1} readily applies to problem \eqref{eq5.6}.  The only 
extra step is to bound the expected value of the  error as follows:
\begin{align}\label{eq6.15}
\bigl\|\mE[u] -\mE_{mc} [U^h_n] \bigr\|_V &\leq  \bigl\|\mE[u] -\mE [U^h_n] \bigr\|_V 
+ \bigl\|\mE[U^h_n] -\mE_{mc} [U^h_n] \bigr\|_V \\
&\leq C \bigl(\rho^{n+1} + h^r \bigr) +   C J^{-\frac12},
\nonumber
\end{align}
where $\mE_{mc}[\cdot] $ denotes the Monte Carlo approximation of the expected value. 
In addition, unlike Approach \#1 above, which is equivalent to the multimodes method 
of \cite{Feng} and imposes the restrictive condition 
$0< \varepsilon <1$ for convergence because $\rho=\epsilon$, we like to show below that using 
Approach \#2 the iterative method also converges for $\varepsilon >1$ because we can have
$\rho< 1$ even $\varepsilon>1$ in this approach.

We reconsider the random diffusion test problem \eqref{boundary_problem2} whose exact solution
is given by \eqref{exact_sol}. It is easy to check that 
\begin{align*}
	\mE[u] = \frac12\Bigl(\frac{1}{\varepsilon} - \frac{1}{\varepsilon^2}\ln(1+\varepsilon)\Bigr)(x-x^2).
\end{align*}
The parameter-dependent problem \eqref{eq5.6} now becomes
\begin{align*}
-\frac{d}{dx}\Bigl(\bigl(1+\varepsilon X(\omega_j)\bigr)\frac{d u(\omega_j,x)}{dx}\Bigr) &= X(\omega_j), \qquad 0 < x < 1,\\\nonumber
u(\omega_j,0) =0,\quad u(\omega_j,1) &= 0,
\end{align*}
for $j = 1, \cdots, J$.

To fit the setup of Section \ref{sec-6.1}, we also define the following bilinear forms:
\begin{alignat*}{2}
&A(\omega_j; u,v) = \bigl(a(\omega_j,\cdot) u',  v'\bigr), &&\qquad F(\omega_j;v) = \bigl(X(\omega_j),v\bigr),\\
&A_0(u,v) = \bigl( a_0\, u', v'\bigr), &&\qquad 
A_1(\omega_j; u,v) = \bigl(\eta(\omega_j,\cdot) u', v'\bigr),
\end{alignat*}
where $a(\omega_j, x) = 1 + \varepsilon X(\omega_j)$; $a_0$ can be chosen as in Section \ref{sec-6.1}. For example, we can take $\displaystyle a_0(x) = \max_{1\leq j\leq J} \,a(\omega_j,x)$ or $a_0(x) = \mE_{mc}[a(\cdot,x)]$. Then, $\eta(\omega_j,x) = a(\omega_j,x) - a_0(x)$ will automatically satisfy the relative dominant condition (ii) with 
\begin{align*}
\displaystyle	\rho = \frac{\max\limits_{1\leq j \leq J}\max\limits_{x\in D}|\eta(\omega_j,x)|}{\min\limits_{x\in D}\,a_0(x)} < 1.
\end{align*}
We also introduce the following error functions:
\begin{align*}
	\mathcal{E}_{L^2}^n = \|\mE[u] - \mE_{mc}[U_n^h]\|_{L^2},\quad  \mathcal{E}_{H^1}^n = \|\mE[u] - \mE_{mc}[U_n^h]\|_{H^1}.
\end{align*}

We now want to verify the convergence rate given in \eqref{eq6.15} for the iterative method.
Notice that there are three terms in the error estimate \eqref{eq6.15}. Since the 
error (i.e., the term $C J^{-\frac12}$) due to the Monte Carlo method is standard and we omit it here. In order to neglect this Monte Carlo error, we choose a large sample number $J =10^6$ in all the numerical tests below. 

To verify the finite element method error term $O(h)$, we consider $\varepsilon = 2.0$ 
and fix $n =10$ then choose different values of the mesh size $h$. It is easy to check that 
$\rho = 0.5003$ when $a_0(x) = \mE_{mc}[a(\cdot,x)]$ and $\rho = 0.6667$ when $ a_0(x) = \max_{1\leq j\leq J} a(\omega_j,x)$. In both cases, $\rho<1$, which satisfies the relative dominant assumption in the convergence theorem. 
Table \ref{tab6.8} displays the $H^1$- and $L^2$-norm errors for the first choice of $a_0$ and  Table \ref{tab6.9} shows the errors for the second choice of $a_0$. We observe a convergence rate $O(h)$ for the $H^1$-norm errors in both cases which is consistent with our error estimate in \eqref{eq6.15}. In addition, we observe that the $L^2$-norm errors exhibit an $O(h^2)$ order
of the convergence. 
\begin{table}[htp]\small
	\renewcommand{\arraystretch}{1.4}
	\centering
	\begin{tabular}{ |l|l|l|l|l| } 
		\hline
		$h$      & $\mathcal{E}_{H^1}^n$ & Order  & $\mathcal{E}_{L^2}^{n}$  & Order   \\ 
		\hline
		$0.2$    &$1.17\times 10^{-2}$  &  & $7.2850 \times 10^{-4}$ &   \\ 
		$0.1$    & $6.20\times 10^{-3}$ &$0.93$  &$1.8724\times 10^{-4}$  & $1.96$  \\ 
		$0.05$    &$3.20\times 10^{-3}$  &$0.96$  & $4.2524\times 10^{-5}$  &$2.14$  \\ 
		$0.025$ & $1.60\times 10^{-3}$ &$1.00$  &$6.9804\times 10^{-6}$  &$2.61$   \\ 
		\hline
	\end{tabular}
	\caption{Approach \#2: Errors $\mathcal{E}_{H^1}^n$ and $\mathcal{E}_{L^2}^n$, for $\varepsilon = 2.0$, $a_0 = \mE_{mc}[a] = 2.0011$, $\rho = 0.5003$, $n =10$ and different $h$.}
	\label{tab6.8}
\end{table}
\begin{table}[htp]\small
	\renewcommand{\arraystretch}{1.4}
	\centering
	\begin{tabular}{ |l|l|l|l|l| } 
		\hline
		$h$      & $\mathcal{E}_{H^1}^n$ & Order  & $\mathcal{E}_{L^2}^{n}$  & Order   \\ 
		\hline
		$0.2$    &$1.17\times 10^{-2}$  &  & $7.3323 \times 10^{-4}$ &   \\ 
		$0.1$    & $6.20\times 10^{-3}$ &$0.93$  &$1.9227\times 10^{-4}$  & $1.93$  \\ 
		$0.05$    &$3.20\times 10^{-3}$  &$0.96$  & $4.7262\times 10^{-5}$  &$2.02$  \\ 
		$0.025$ & $1.60\times 10^{-3}$ &$1.00$  &$9.9898\times 10^{-6}$  &$2.24$   \\ 
		\hline
	\end{tabular}
	\caption{Approach \#2: Errors $\mathcal{E}_{H^1}^n$ and $\mathcal{E}_{L^2}^n$, for $\varepsilon = 2.0$, $\displaystyle a_0(x) = \max_{1\leq j\leq J} a(\omega_j,x) = 3.0$, $\rho = 0.6667$, $n =10$ and different $h$.}
	\label{tab6.9}
\end{table}

\begin{figure}[h]
	\begin{center}
		\includegraphics[scale=0.26]{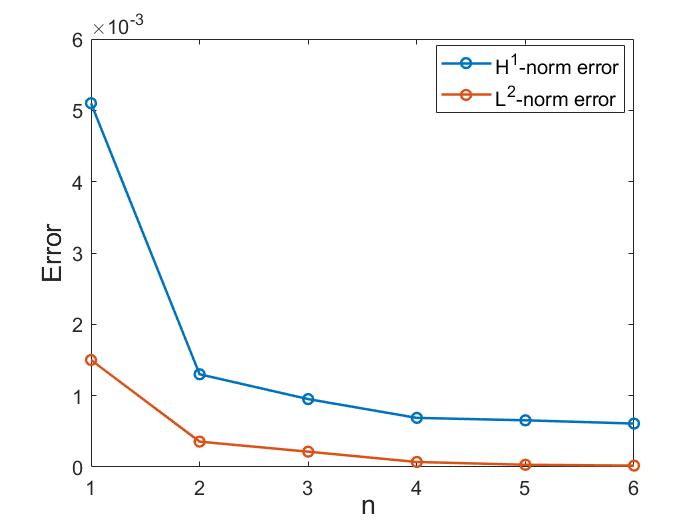}
		\includegraphics[scale=0.26]{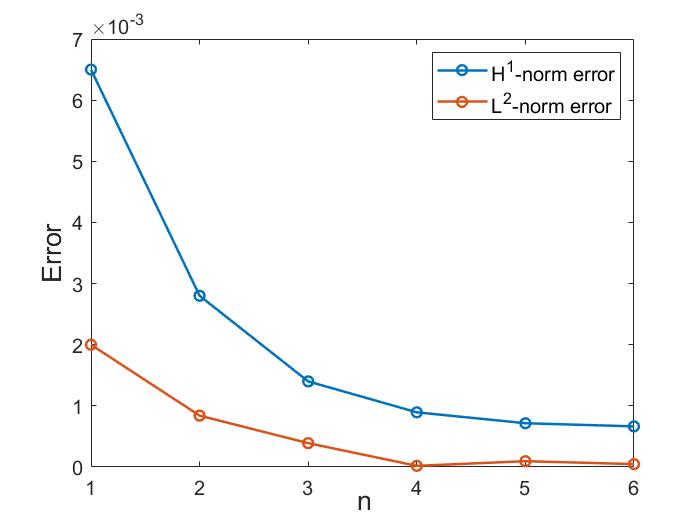}
		\caption{Approach \#2: Errors $\mathcal{E}_{H^1}^n$ and $\mathcal{E}_{L^2}^n$ with $a_0(x) = \mE_{mc}[a(\cdot, x)] = 2.0011$ (left) and $\displaystyle a_0(x) = \max_{1\leq j\leq J} a(\omega_j,x) = 3.0$ (right) for $h = 0.01, \varepsilon = 2.0$, ${\tt tol} = 10^{-4}$ and various $n$.}\label{H1_L2_2ndapp}
	\end{center}
\end{figure}

To verify the iteration error term $O(\rho^{n+1})$, we fix $h = 0.01$, $\varepsilon = 2.0$ and 
also consider two cases of $a_0$ as in the previous test. A stopping criteria   $\|\mE_{mc}\bigl[U^h_{n+1} - U^h_n\bigr]\|_{H^1} < {\tt tol} := 10^{-4}$ is used to terminate
the iteration. 
Figure \ref{H1_L2_2ndapp} displays the $H^1$- and $L^2$-norm errors 
of the computed solutions. We observe that the errors are decreasing as the number of iterations increases. Figure \ref{H1_L2_2ndapp} also shows that only $6$ iterations are needed 
to trigger the stopping criterion.  

\begin{figure}[h]
	\begin{center}
		\includegraphics[scale=0.26]{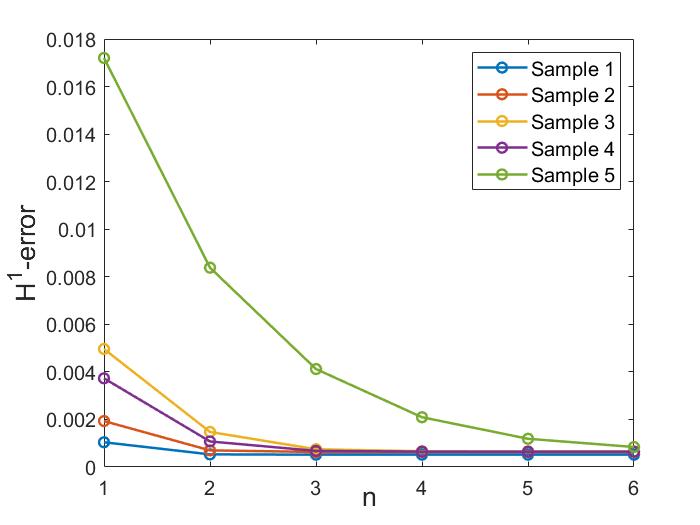}
		\includegraphics[scale=0.26]{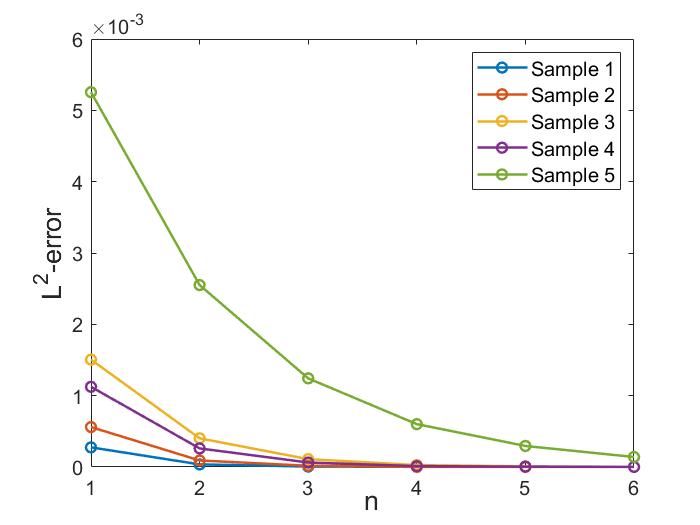}
		\caption{Approach \#2: Errors $\mathcal{E}_{H^1}^n$ (left) and $\mathcal{E}_{L^2}^n$ (right) decays for first $5$ samples with $ \varepsilon =2.0, a_0(x) = \mE_{mc}[a(\cdot,x)] = 2.001, \rho = 0.5003, h =0.01$ and various $n$.}\label{5samples_2ndapp}
	\end{center}
\end{figure}

Finally, we plot the errors of the first five Monte Carlo samples in Figure \ref{5samples_2ndapp} for $\varepsilon = 2.0, a_0(x) = \mE_{mc}[a(\cdot,x)], h = 0.01$ and various $n$. 
As expected,  the errors are different for different samples and they all are also different from 
the expected value in Figure \ref{H1_L2_2ndapp}, but they all are decreasing as the number of iterations $n$ increases.

\subsection{1-D random convection-diffusion problems}\label{sec-6.3}
In this subsection, we consider the following 1-D randomized double-glazing problem \cite{IFISS}:
\begin{align}\label{eq6.16} 
	-\frac{d}{dx}\Bigl(\bigl(1+ \varepsilon X(\omega)\bigr)\frac{d u(\omega,x)}{dx}\Bigr) &+ 100\bigl(1+\varepsilon X(\omega)\bigr) \frac{d u(\omega, x)}{d x} = f(\omega,x),\\\nonumber
	u(\omega,0) &= 0,\qquad u(\omega,1) = 0,
\end{align} 
where $0\leq x\leq 1, \varepsilon > 0$ and $X(\omega)$ is a uniformly distributed random variable on a probability space $(\Ome, \mathcal{B},\mP)$ with $\mathcal{B}$ is the Borel $\sigma$-algebra. In addition, we choose $f(\omega,x) = X(\omega)(51-100x)$ so that the exact solution is again given by
\begin{align}
	u(\omega,x) = \frac{X(\omega)}{2(1+\varepsilon X(\omega))}(x-x^2).
\end{align}
Its expected value is already given in previous subsection. 

To apply the iterative method, we adopt Approach \#2 proposed in Section \ref{sec-6.1} to \eqref{eq6.16}. Namely, we consider the following parameter-dependent problem:
\begin{align}\label{eq6.17}
	-\frac{d}{dx}\Bigl(\bigl(1+ \varepsilon X(\omega_j)\bigr)\frac{d u(\omega_j,x)}{dx}\Bigr) &+ 100\bigl(1+\varepsilon X(\omega_j)\bigr) \frac{d u(\omega_j, x)}{d x} = f(\omega_j,x),\\\nonumber
	u(\omega_j,0) &= 0,\qquad u(\omega_j,1) = 0,
\end{align} 
for $j = 1,\cdots, J$. To fit the setup of Section \ref{sec-6.1}, we introduce the following coefficient functions:
\begin{alignat*}{3}
	&a(\omega_j,x) = 1 + \varepsilon X(\omega_j),&&\qquad b(\omega_j,x) = 100(1 + \varepsilon X(\omega_j))\\\nonumber
	&a_0(x) = \mE_{mc} [a(\cdot,x)],&&\qquad b_0(x) = \mE_{mc}[b(\cdot,x)],\\\nonumber
	&\eta_a(\omega_j,x) = a(\omega_j,x) - a_0(x),&&\qquad \eta_b(\omega_j,x) = b(\omega_j,x) - b_0(x),
\end{alignat*}
which result in the following bilinear forms:
\begin{align*}
	&A(\omega_j; u,v) = \bigl(a(\omega_j,\cdot) u',  v'\bigr) + \bigl(b(\omega_j,\cdot)u',v\bigr), \qquad F(\omega_j;v) = \bigl(f(\omega_j,\cdot),v\bigr),\\
	&A_0(u,v) = \bigl( a_0\, u', v'\bigr) + \bigl(b_0u',v\bigr), \quad 
	A_1(\omega_j; u,v) = \bigl(\eta_a(\omega_j,\cdot) u', v'\bigr) + \bigl(\eta_b(\omega_j,\cdot) u', v\bigr).
\end{align*}
Then, the relative dominant number $\rho$ in the convergence criteria of Section \ref{sec-3}
can be estimated by
\begin{align*}
	\displaystyle	\rho \approx \frac{\max\limits_{1\leq j \leq J}\max\limits_{x\in D}(|\eta_a(\omega_j,x)|+|\eta_b(\omega_j,x)|)}{\min\limits_{x\in D}\,a_0(x)}.
\end{align*}

In this test, we use $J = 10^4$ number of samples for the Monte Carlo method, the mesh size $h = 0.01$ and the maximum number of iterations $n = 10$. We also select $\varepsilon = 0.2, 0.005$ in \eqref{eq6.17}.   
\begin{figure}[h]
	\begin{center}
		\includegraphics[scale=0.26]{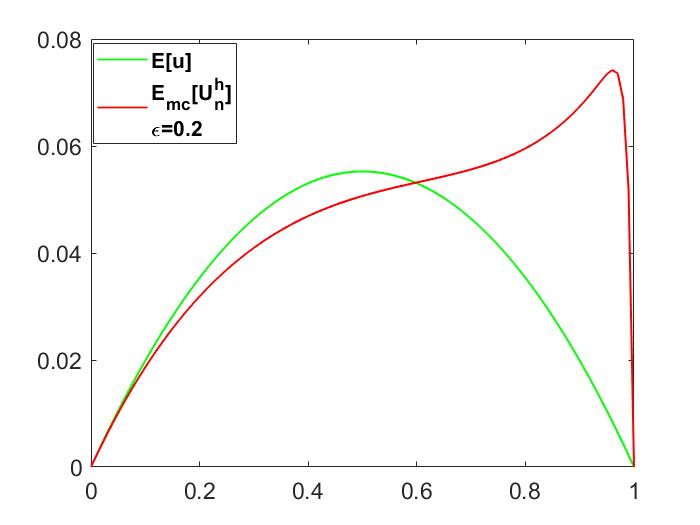}
		\includegraphics[scale=0.26]{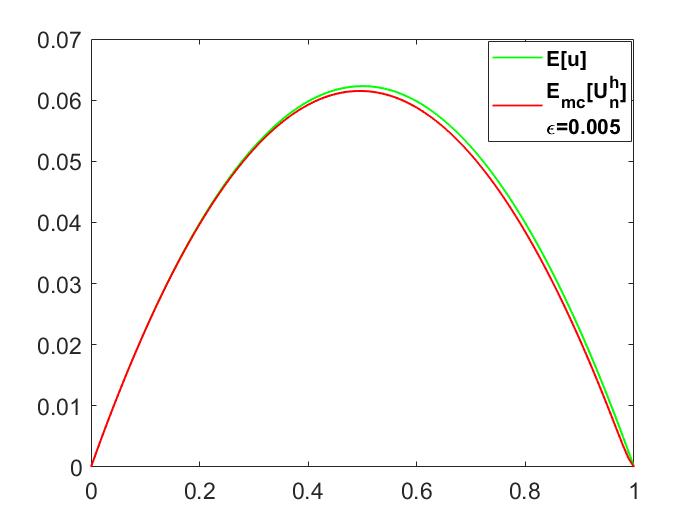}
		\caption{1-D random convection-diffusion problem: Plots of $\mE[u]$ and $\mE_{mc}[U_n^h]$ with $h = 0.01, n =10$ and $\varepsilon = 0.2$ (left) and $\varepsilon = 0.005$(right).}\label{convec_diff_solution}
	\end{center}
\end{figure}

\begin{table}[ht]\small
	\renewcommand{\arraystretch}{1.4}
	\centering
	\begin{tabular}{ |l|l|l|l| } 
		\hline
		$\varepsilon$      & $\rho$ & $\mathcal{E}^n_{H^1} = \|\mE[u] - \mE_{mc}[U_n^h]\|_{H^1}$ &  $\mathcal{E}^n_{L^2} = \|\mE[u] - \mE_{mc}[U_n^h]\|_{L^2}$\\
		\hline
		$0.2$& $9.1912$ & $0.5328$&$2.2100\times 10^{-2}$\\
		\hline
		$0.005$&$0.2522$&$0.0117$&$9.6119\times 10^{-4}$\\
		\hline
	\end{tabular}
	\caption{1-D random convection-diffusion problem: Errors $\mathcal{E}_{H^1}^n$ and $\mathcal{E}_{L^2}^n$ with $h = 0.01, n = 10$ and various $\varepsilon$.}
	\label{tab6.10}
\end{table}

Figure \ref{convec_diff_solution} shows the expected values $\mE_{mc}[U_n^h]$ and $\mE[u]$ of the computed and exact solutions and Table \ref{tab6.10} displays the expected values of 
the $H^1$- and $L^2$-norm errors. 
We observe from Figure \ref{convec_diff_solution} and Table \ref{tab6.10} that the proposed iterative method performs well for problem \eqref{eq6.16} when $\varepsilon$ is sufficiently small, but the errors become larger as $\varepsilon$ increases, which is expected. In addition, due to the contribution of the convection term in \eqref{eq6.16}, the relative dominance parameter $\rho>1$ when $\varepsilon = 0.2$, which explains the
poor performance of the iterative method in this case.

	\subsection{2-D randomized double-glazing problems}\label{sec-6.4}
In this subsection, we consider the random perturbation of the double-glazing problem 
given in \cite{IFISS}. Specifically, we consider
\begin{alignat}{2}\label{double_glazing}
	-\delta \Delta u + \bb\cdot\nab u &= f(\omega,x,y) &&\qquad \mbox{ in } D = (0,1)^2\\
	u(\omega,\cdot) &= 0 &&\qquad \mbox{ on } \p D,
\end{alignat}
where $\bb(\omega,x,y) := (1+\varepsilon X(\omega))\bigl(2y(1-x^2), -2x(1-y^2)\bigr)$ and $0< \delta << |\bb|$, $\varepsilon >0$. $X$ is a uniformly distributed random variable on $[0,1]$. It is easy to check that $\div \bb = 0$.  We set the right hand-side force function
\begin{align*}
f(\omega;x,y) &= \frac{\delta X(\omega)}{1 + \varepsilon X(\omega)}(y - y^2 + x - x^2) + X(\omega)y(y-y^2)(1-x^2)(1-2x)\\
&\qquad  - X(\omega)x(x-x^2)(1-y^2)(1-2y)
\end{align*}
so that the exact solution is  given by
\begin{align*}
	u(\omega;x,y) = \frac{X(\omega)}{2(1+\varepsilon X(\omega))}(x-x^2)(y-y^2).
\end{align*} 
Thus,
\begin{align*}
	\mE[u] = \frac12\Bigl(\frac{1}{\varepsilon} - \frac{1}{\varepsilon^2}\ln(1+\varepsilon)\Bigr)(x-x^2)(y-y^2).
\end{align*}
 
We again adopt Approach \#2 of Section \ref{sec-6.1} to solve \eqref{double_glazing} with $\delta = 0.1$ using 
the proposed iterative method with $J = 10^4$ number of Monte Carlo samples. Consequently, we need to solve the following parameter-dependent problem: for $1\leq j \leq J$
\begin{alignat}{2}
	-\delta \Delta u(\omega_j,x,y) + \bb(\omega_j,x,y)\cdot\nab u(\omega_j,x,y) &= f(\omega_j,x,y) &&\quad \mbox{ in } D = (0,1)^2\\
	u(\omega_j,\cdot) &= 0 &&\quad\mbox{ on } \p D. \nonumber
\end{alignat}
 
To fit the setup of Section \ref{sec-6.1}, we define
\begin{align*}
	\bb_0(x,y) &:= \mE_{mc}[\bb(\cdot,x,y)]=2.001 \bigl(2y(1-x^2), -2x(1-y^2)\bigr),\\
	\pmb{\eta}(\omega_j,x,y) &:= \bb(\omega_j,x,y) - \bb_0(x,y),
\end{align*}
and the following bilinear forms and functional:
\begin{alignat*}{2}
	&A(\omega_j; u,v) = \delta\bigl(\nab u,  \nab v\bigr) + \bigl(\bb\cdot\nab u,v\bigr), &&\quad F(\omega_j;v) = \bigl(f,v\bigr),\\
	&A_0(u,v) = \delta\bigl( \nab\, u, \nab v\bigr) + \bigl(\bb_0\cdot\nab u,v\bigr), &&\quad
	A_1(\omega_j; u,v) = \bigl(\pmb{\eta}\cdot \nab u, v\bigr).
\end{alignat*}
 
\begin{table}[htp]\small
	\renewcommand{\arraystretch}{1.4}
	\centering
	\begin{tabular}{ |l|l|l|l|l| } 
		\hline
		$h$      & $\mathcal{E}_{H^1}^n$ & Order  & $\mathcal{E}_{L^2}^{n}$  & Order   \\ 
		\hline
		$0.2$    &$5.4020\times 10^{-3}$  &  & $3.4360 \times 10^{-4}$ &   \\ 
		$0.1$    & $2.7339\times 10^{-3}$ &$0.9825$  &$8.2969\times 10^{-5}$  & $2.0501$  \\ 
		$0.05$    &$1.3709\times 10^{-3}$  &$0.9958$  & $1.7666\times 10^{-5}$  &$2.2316$  \\ 
		$0.025$ & $6.8635\times 10^{-4}$ &$0.9806$  &$3.8150\times 10^{-6}$  &$2.2112$   \\ 
		\hline
	\end{tabular}
	\caption{2-D random convection-diffusion problem: Errors $\mathcal{E}_{H^1}^n$ and $\mathcal{E}_{L^2}^n$ with $\varepsilon = 2.0$,  
		$n =10$ and various $h$.}\label{table6.10}
\end{table}

\begin{figure}[h]
	\begin{center}
		\includegraphics[scale=0.26]{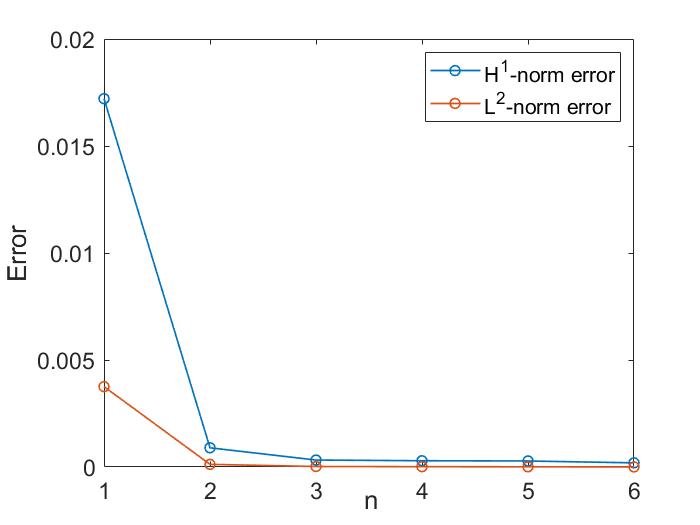}
		\caption{2-D convection-diffusion problem: Errors $\mathcal{E}_{H^1}^n$ and $\mathcal{E}_{L^2}^n$ with  
			$h = 0.01, \varepsilon = 2.0$, ${\tt tol} = 10^{-4}$ and various $n$.}\label{fig6.14}
	\end{center}
\end{figure}

Table \ref{table6.10} and Figure \ref{fig6.14} show respectively the convergence orders 
of the spatial 
discretization error and the iteration errors in both $H^1$- and $L^2$-norm. 
We observe that the proposed iterative method performs well for the random convection-dominated convection-diffusion problem \eqref{double_glazing}. 
From Table \ref{table6.10} we see that the optimal convergence orders are achieved 
for the linear finite element method. In addition, Figure \ref{fig6.14} shows 
that the $H^1$- and $L^2$-norm errors of the iterative solution rapidly decrease 
as the number of iterations increases.





\end{document}